\documentclass[final]{siamart190516}
 \usepackage{marginfix}
 \usepackage{amsfonts}

 \usepackage{amsmath}
 \usepackage{amssymb}
 \usepackage{enumerate}
 \usepackage{graphicx}  
 \usepackage{color}
 \usepackage{longtable}
 \usepackage{multirow}
 \usepackage{float}
 \usepackage{booktabs}
 \usepackage{makecell}
 \usepackage{supertabular}
 \usepackage{algpseudocode}
 \usepackage{algorithmicx,algorithm}
 \usepackage{caption}
 \usepackage{subfigure}
 \usepackage{hyperref}
 
 \numberwithin{equation}{section}
 \newtheorem{remark}{Remark}[section]

 \allowdisplaybreaks[4]

\newcommand{\cT}{\mathcal{T}}

\begin{document}

\title{A fast multilevel dimension iteration algorithm for high dimensional numerical integration}

\headers{A FAST ALGORITHM FOR HIGH DIMENSIONAL NUMERICAL INTEGRATIPON}{HUICONG ZHONG AND XIAOBING FENG} 

\author{
	Xiaobing Feng\thanks{Department of Mathematics, The University of Tennessee,  Knoxville, TN, 37996 (\email{xfeng@utk.edu}}. 
	\and
	Huicong Zhong\thanks{School of Mathematics and Statistics, Northwestern Polytechnical University, Xi’an, Shaanxi 710129, China (\email{huicongzhong@mail.nwpu.edu.cn})}.
}

\maketitle

\begin{abstract}
	In this paper we propose and study a fast multilevel dimension iteration (MDI) algorithm for computing arbitrary $d$-dimensional integrals based on the tensor product approximations. It reduces the computational complexity (in terms of the CPU time) of a tensor product method from the exponential order $O(N^d)$ to the polynomial order  {\color{black} $O(d^3N^2)$ or better}, where $N$ stands for the number of quadrature points in each coordinate direction. As a result, the proposed MDI algorithm effectively  circumvents the curse of the dimensionality of  tensor product methods for high dimensional numerical integration. The main idea of the proposed MDI algorithm is to compute the function evaluations at all integration points in cluster and iteratively along each coordinate direction, so lots of computations   
	for function evaluations can be reused in each iteration. This idea is also applicable to any quadrature rule whose integration points have a lattice-like structure.
	
\end{abstract}

\begin{keywords}
	 Multilevel dimension iteration (MDI), high dimensional integration, numerical quadrature rules, tensor product methods,  Monte Carlo methods. 
\end{keywords}

\begin{AMS}
65D30, 65D40, 65C05, 65N99
\end{AMS}

\section{Introduction}\label{sec-1}
Numerical integration is one of the most fundamental building blocks in computational mathematics and in computational science at large. Many numerical methods (or quadrature rules) had been well developed as documented in numerical analysis textbooks (cf. \cite{Burden-Faires, Stoer-Bulirsch} and the references therein). 
They are tensor-product-type methods and work very well for computing integration in  low dimensions. However, they all become prohibitively expensive in the high-dimensional case because the number of
function evaluations grows as $O(N^d)$ (hence, the CPU time for computing them is expected to grow even faster),
where $d$ and $N$ denote respectively the dimension of the integration domain 
and the number of the quadrature points in each coordinate direction. This exponential
growth of the computational complexity is known as ``the curse of the dimensionality"
(cf. \cite{BG14, DKS13}). 

To circumvent or lessen ``the curse of the dimensionality", various improved numerical integration methods have been developed in the literature.  Among them are 
the Monte Carlo (MC) method and its variants including quasi- and multilevel Monte Carlo (QMC) methods  \cite{Caflisch98,DKS13,HMNR10,KSS11,Ogata89,DKS13,Lu04,Wipf13}, and sparse grid (SG) methods \cite{BG14,GH10,Dos12}, and deep neural network-based methods \cite{E_Yu18, HJE18, LMMK21, SS18,Xu20}.  Although those improved methods are fundamentally different in their 
design ideas and mechanisms, they share a common strategy  that is 
to reduce the number function evaluations compared to the tensor-product methods while maintaining the reasonable degree of accuracy. As expected, such a strategy 
works to some extent for medium dimensions (i.e., $d\approx 100$), but may not work for very high dimensions (i.e., $d\approx 1000$) because the required number of function evaluations still grows very rapidly  for large $d$. The situation is even direr if one wants to solve partial differential equations (PDEs) in high dimensions. Consequently, developing faster and efficient numerical integration methods remains critical for tackling more challenging problems
arising from mathematical finance, image processing, economics and data science.

The primary goal of this paper is to develop a fast algorithm, called 
{\em multilevel dimension iteration (MDI)}, for high dimensional numerical 
integration.  Unlike the MC, QMC, SG and deep neural network (DNN) methods, 
the proposed MDI algorithm is not aiming to providing a new numerical integration method per se, 
instead, it is an acceleration algorithm for an efficient implementation of any  tensor-product-type method. Thus, the MDI is not a ``discretization" method but a ``solver" 
(borrowing the numerical PDE terminologies). A well suited analogy would be high order 
polynomial evaluations, that is, to compute 
$p_0:=p(x_0)=a_kx_0^k+a_{n-1}x_0^{n-1}+\cdots + a_1x_0+a_0$ for a given real number
input $x_0$. It is well known that such a high order polynomial evaluation on a computer is notoriously unstable, inaccurate (due to roundoff errors) and expensive, however, those difficulties can be easily overcome 
by a simple nested iteration (or Horner's algorithm. cf. \cite{Burden-Faires,Stoer-Bulirsch}), namely, set $p_0:=a_k$ and for $j=k, k-1, \cdots, 1$, compute $p_0:=p_0x_0+a_{j-1}$.
From the cost saving and efficiency point view, the reason for the nested iteration to be efficient and fast is that it reuses many multiplications involving $x_0$ compared to the direct evaluations of each term in $p(x_0)$. Conceptually, this is exactly the approach adopted by the proposed MDI algorithm, i.e., to reuse computations of 
the function evaluations in a tensor product method as much as possible 
to save the computation cost and hence to make it efficient and fast.  
A key observation is that the function evaluations of every tensor product method (including SG) involve a lot of computation in each coordinate direction which 
can be shared because each coordinate $\xi_j$ of every integration point $\xi=(\xi_1,\xi_2,\cdots, \xi_d) \in \mathbb{R}^d$ is shared by many other 
integration points due to their tensor product structure. This observation motivates us
to compute the required function evaluations in cluster and iteratively in each coordinate
direction instead
of computing them at the integration points independently, which is 
exactly the key idea of the proposed MDI algorithm. In other words, our MDI algorithm is based on 
a numerical Fubini's approach to compute the summation (and function evaluations) defined by a given tensor product method.

The remainder of the paper is organized as follows. In Section \ref{sec-2} we 
introduce our MDI algorithm first in two and three dimensions to explain the main   
ideas of the algorithm, and then generalize them to arbitrary dimensions. In Section \ref{sec-3} we present various numerical experiments to test the performance of
the proposed MDI algorithm and to do various performance comparisons with the classical MC method. It shows that the MDI algorithm (regardless the choice of the underlying tensor product method) is 
faster than the classical MC method in low and medium dimensions (i.e., $d\approx 100$), much faster in 
very high dimensions (i.e., $d\approx 1000$), and succeeds even when the MC method fails. 
In Section \ref{sec-4} we provide numerical experiments to gauge the 
influence of parameters in the proposed MDI algorithm, including the dependence on
 choices of the underlying tensor product method and the iteration step size. 
In Section \ref{sec-5}, we use the computation techniques to find out the computational 
complexity of the MDI algorithm. This is done by discovering the relationship between CPU time and  
dimension by the standard regression technique. 
it shows that the CPU time grows in the polynomial order $O(d^3N^2)$ at most.
Furthermore,  numerical experiments are designed to test the limit of the proposed MDI algorithm, it can compute integrals  on standard desktop computers 
in medium dimensions easily and in very high dimensions quickly. Finally, we complete the paper with a few concluding remarks given in Section \ref{sec-6}.

\section{Methodology: the MDI algorithm}\label{sec-2}
Let $\Omega\subset \mathbb{R}^d$ ($d\geq 1$) be a bounded rectangular domain. Without loss of the generality, unless mentioned otherwise, we assume $\Omega=[0,1]^d$.  
Let $g:\overline{\Omega}\to \mathbb{R}$ denote a generic continuous function on 
$\overline{\Omega}$ (hence, $g$ has a pointwise value at every $x=(x_1,x_2,\cdots,x_d)\in \overline{\Omega}$). Then the central issue to be addressed in this paper is to evaluate 
\begin{equation}\label{eq1}
I(g,\Omega):=	\int_{\Omega} g(x)\, dx 
\end{equation}
accurately and efficiently for $d>>1$. 

\subsection{Tensor product methods}\label{sec-2.1}
In this subsection, we briefly recall the formulation of tensor product methods for approximating \eqref{eq1} and their well-known properties. 

Let $N$ be a positive integer and  $h=\frac{1}{N}$. 
Let $\mathcal{T}_h=\{K\}$ denote the uniform rectangular mesh of $\Omega$ with mesh 
size $h$. Note that each element $K\in\mathcal{T}_h$ is a $d$-rectangle (in fact, a $d$-square of side $h$ when $\Omega$ is a $d$-square domain) and the 
total number of elements in $\mathcal{T}_h$ is $N^d$. Define for $K\in \mathcal{T}_h$
\begin{equation}\label{eq2.1}
I(g,K):=\int_K g(x)\, dx.
\end{equation}
By the summation property of integrals, we have
\begin{equation}\label{eq2.2a}
	I(g,\Omega)= \sum_{K\in \cT_h} I(g,K). 
\end{equation}
Now,  approximating every $I(g,K)$ by a local numerical quadrature rule 
\begin{equation}\label{eq2.2b}
J(g,K) \approx I(g,K) \qquad \forall K\in \cT_h, 
\end{equation}
it then leads to the following global approximation 
\begin{equation}\label{eq2.2c}
J(g,\Omega):=  \sum_{K\in \cT_h} J(g,K) \approx  \sum_{K\in \cT_h} I(g,K) =I(g,\Omega), 
\end{equation}
which is known as the composite method of the chosen local quadrature rule. 

When the local 
rule \eqref{eq2.2b} is constructed by repeatedly using the same one-dimensional quadrature rule in each coordinate direction, then the resulting global rule \eqref{eq2.2c} is called a tensor product quadrature rule for approximating 
$I(g,\Omega)$. 

In this paper we only focus on the following four popular tensor product
quadrature rules: (i) the trapezoidal rule; (ii) the Simpson's rule; (iii) the 
two-point Gaussian rule; (iv) the midpoint rule. It is well known \cite{Burden-Faires,Stoer-Bulirsch} that these four rules have the following error estimates:
{\color{black}
\begin{alignat*}{2}
	   &{\rm (i)}\,\, I(g,\Omega)-J_{\mbox{\tiny{trap}}}(g,\Omega) = O(N^{-2}), 
	   \qquad 	  &&{\rm (ii)}\,\, I(g,\Omega)-J_{\mbox{\tiny{simp}}}(g,\Omega) = O(N^{-4}), \\
	 &{\rm (iii)}\,\, I(g,\Omega)-J_{\mbox{\tiny{gauss}}}(g,\Omega) = O(N^{-4}),  &&
      {\rm (iv)}\,\, I(g,\Omega)-J_{\mbox{\tiny{midp}}}(g,\Omega) = O(N^{-2}). 
\end{alignat*}
}

Clearly, these error estimates hold  only when $g$ is sufficiently smooth. 
Also, the trapezoidal rule is lower order compared to the Simpson and two-point Gaussian rules, however, the midpoint, trapezoidal and Simpson rules are easier to implement and flexible on choosing the integration points. This feature is advantageous in the case when explicit formula of the integrand $g$ is not known. On the other hand, the  two-point Gaussian rule has higher accuracy at the expense of computing the Gaussian points (which is costly when $d>>1$). In Section \ref{sec-4} we shall use numerical experiments to further elaborate this point in the context of the proposed MDI algorithm.

\subsection{Formulation of the MDI algorithm in two dimensions}\label{sec-2.2}
To better understand and to present the idea of the MDI algorithm, 
we first consider the simple 
two dimensional case (i.e., $d=2$ and $\Omega=[0,1]^2$). 
In the two-dimensional case, by Fubini's Theorem we have 
\begin{equation}\label{eq2.6}
	I\bigl(g,[0,1]^2 \bigr) :=\int_{[0,1]^2} g(x)\, dx = \int_{0}^{1} \biggl( \int_{0}^{1} g(x) \,dx_1 \biggr) \,dx_2. 
\end{equation}
It should be noted that the exact evaluation of the above integral (assuming it is doable) by hand 
is often done using Fubini's theorem in calculus. Conceptually, this trivial fact will play an important role in conceiving the idea of our MDI algorithm. 

Let $N\geq 1$ be an integer. Suppose that we have (or choose) the following generic one-dimensional quadrature rule: 
\begin{equation}\label{eq2.6b}
 \int_0^1 \phi(s) \, ds\approx \sum_{i=1}^N w_i \phi(\xi_i) =:J(g,[0,1]),
\end{equation}
where $\{\xi_i\}$ and $\{w_i\}$ denote respectively the nodes and weights of the quadrature rule. 
$\phi$ is a generic (continuous) function on $[0,1]$.  

Then, the corresponding two-dimensional tensor product rule takes the form
\begin{equation} \label{eq2.7}
 J\bigl(g, [0,1]^2\bigr)=\sum_{i,j=1}^{N} w_i w_j g(\xi_{ij}) 
\end{equation}
where the two-dimensional nodes $\xi_{ij}:=(\xi_i, \xi_j)$.   Obviously, the computational complexity of the above quadrature rule is $O(N^2)$

%

Motivated by (and mimicking) the Fubini's formula \eqref{eq2.6}, we rewrite the tensor product rule \eqref{eq2.7} as 
\begin{equation}\label{eq2.8}
J\bigl(g, [0,1]^2\bigr)=\sum_{j=1}^{N} w_j \bigg( \sum_{i=1}^N  w_i g(\xi_{ij}) \biggr) 
=\sum_{j=1}^N w_j g_1(\xi_j),
\end{equation}
where  
\begin{equation}\label{eq2.8b}
	g_1(s) :=\sum_{i=1}^{N} w_i g\bigl((\xi_i,s) \bigr).
\end{equation}
We note that the evaluation of $g_1(\xi_j)$ is amount to applying the 1-d formula \eqref{eq2.6b} 
to approximate the integral $\int_0^1 g\bigl((x_1, \xi_j))\, dx_1.$ 
However, the values of $\{g_1(\xi_j) \}$ will not be computed by the 1-d quadrature rule
in our MDI algorithm, instead, $g_1$ is formed as a symbolic function, so the 1-d 
quadrature rule can be called again.

Evidently, \eqref{eq2.8} is a discrete analogue of the Fubini's formula \eqref{eq2.6b}, 
hence, we refer \eqref{eq2.8} as a {\em discrete Fubini's formula} in the rest of this paper. 
This simple formula has a significant computational benefit because it does all computations 
which involve the first (i.e., $x_1$) components of all two-dimensional integration nodes/points 
first and save them in terms of the symbolic function $g_1$ defined by \eqref{eq2.8b}, 
then the final function evaluations in the tensor product formula are done by evaluating $g_1$ at the second (i.e., $x_2$) component of all two-dimensional integration nodes/points.

Let $W$ and $X$ denote the weight and node vectors of a selected 1-d numerical quadrature rule 
on $[0,1]$ and we use a parameter $r$ to indicate the quadrature rule. The following algorithm implements the above discrete Fubini's formula. 










\begin{algorithm}[H]
		\caption{2d-MDI($g$, $\Omega, N, r$)} 
		\label{alg:1}
		\hspace*{0.02in} {\bf Inputs:} 
	    $g, \Omega, N, r.$ \\
		\hspace*{0.02in} {\bf Output:}
	    $J=J(g,\Omega)$.
		\begin{algorithmic}[1]
			\State Initialize $g_1=0$, $J=0$.
	     	\For{$i=1:N$} 	 
			\State $g_1=g_1+W(i) g((X(i), \cdot))$.
			\EndFor	
			\For{$j= 1:N$} 
			\State $J=J+W(j) g_1(X(j))$.
			\EndFor\\		
		\Return $J$.
		\end{algorithmic}
\end{algorithm}
We note that the first do-loop forms the symbolic function $g_1$ which encodes all  
computations involving the $x_1$-components of all integration points. The second do-loop 
evaluates the 1-d quadrature rule for the function $g_1$.  
As mentioned above, in this paper we only focus on the four well-known 1-d 
quadrature rules: (i) the trapezoidal rule; (ii) the Simpson's rule; (iii) the 
two-point Gaussian rule; (iv) the midpoint rule. They will be represented respectively by $r=1,2,3,4$.

\subsection{Formulation of the MDI algorithm in three dimensions}\label{sec-2.3}
In the subsection we shall formulate the MDI algorithm in the 3-d case. 
Since the main idea is similar to that of the 2-d case, we shall only highlight 
its main steps. 

Applying the 1-d quadrature rule \eqref{eq2.6b} in each of three coordinate directions, we
readily obtain the following 3-d tensor product rule for approximating integral $I(g, [0,1]^3)$:  
\begin{equation} \label{eq2.9}
J\bigl(g, [0,1]^3\bigr)=\sum_{i,j,k=1}^{N} w_i w_j w_k g(\xi_{ijk}) 
\end{equation}
where the three-dimensional nodes $\xi_{ijk}:=(\xi_i, \xi_j,\xi_k)$.   Obviously, the computational complexity of the above formula is $O(N^3)$

Again, by Fubini's Theorem we have 
\begin{equation}\label{eq2.10}
I\bigl(g,[0,1]^3 \bigr) =\int_{[0,1]^3} g(x)\, dx =    \int_{[0,1]^2} \biggl(\int_{0}^{1} g(x)\, dx_1 \biggr)\,dx'  , 
\end{equation}
where $x'=(x_2,x_3)$. 
Mimicking the above Fubini's formula, we rewrite \eqref{eq2.9} as 
\begin{equation}\label{eq7}
	J\bigl(g, [0,1]^3\bigr)=   \sum_{j,k=1}^N w_jw_k \biggl( \sum_{i=1}^{N}  w_i g(\xi_{ijk}) \biggr) 
	= \sum_{j,k=1}^{N} w_j w_k g_2(\xi_j,\xi_k),
\end{equation}
where 
\begin{equation}\label{eq7b}
	g_2(s,t) :=\sum_{i=1}^{N} w_i g((\xi_i,s,t)).
\end{equation}
Once again, it should be noted that $g_2$ will be formed as a symbolic function in our MDI algorithm and the right-hand side of \eqref{eq7} is viewed as a 2-d tensor product formula for $g_2$,  
it  can be computed either directly or recursively by 
using \textbf{Algorithm 2.1}. Below we present our MDI algorithm for 
implementing the recursive strategy.  


 





\begin{algorithm}[H]
	\caption{3d-MDI($g$, $\Omega, N, r$)} 
	\label{alg:2}
	\hspace*{0.02in} {\bf Inputs:} 
	$g, \Omega, N, r.$ \\
	\hspace*{0.02in} {\bf Output:}
	$J=J(g,\Omega)$.
	\begin{algorithmic}[1]
		\State Initialize $g_2=0$, $J=0$.
		\For{$i=1:N$} 	 
		\State $g_2=g_2+W(i) g((X(i), \cdot, \cdot))$.
		\EndFor
		\State $\Omega_2=P_3^2 \Omega$.
		\State $J=$2d-MDI$(g_2, \Omega_2, N, r)$. \\
		\Return $J$.
	\end{algorithmic}
\end{algorithm}
\noindent
where $P_3^2$ denotes the orthogonal projection (or natural 
embedding): $x=(x_1,x_2,x_3)$ $\to x'=(x_2,x_3)$,
$W$ and $X$ stand for the weight and node vectors of the selected 1-d quadrature rule. 

From \textbf{Algorithm 2.2} we already can see the procedure of the MDI algorithm. 
It is based on the two main ideas: (i) to use the discrete 
Fubini's formula to reduce the computation of the tensor product sum into the computation 
of a lower dimensional tensor product sums, which allow us to call recursively a lower dimensional MDI algorithm;
(ii) the function evaluations are done in cluster in each coordinate direction during the 
dimension iteration/reduction, which is the main reason for a significant computational saving
due to reusing lots of computations, compared to the standard pointwise function evaluations which treat 
all the integration points independently and do not reuse any computation.

\subsection{Formulation of the MDI algorithm in arbitrary d-dimensions}\label{sec-2.4}
The goal of this subsection is to extend the 2- and 3-d MDI algorithms to arbitrary d-dimensions. To the end, we first recall a more general version of Fubini's Theorem stated as follows:
\begin{equation}\label{eq2.13}
I(g, \Omega)=\int_{\Omega} g(x)\, dx = \int_{\Omega_{d-m}} \biggl(  \int_{\Omega_m} g(x)\, dx'' \biggr)
\, dx',
\end{equation}
where $1\leq m<d$, $\Omega=[0,1]^d, \Omega_m=Q_d^m \Omega=[0,1]^m$ and $\Omega_{d-m}=P_d^{d-m} \Omega =[0,1]^{d-m}$ in which 
$Q_d^m$ and $P_d^{d-m}$ denote respectively the orthogonal projections 
(or natural embeddings): $x=(x_1,x_2, \cdots, x_d)\to x''=(x_1,x_2, \cdots, x_m)$ and $x=(x_1,x_2, \cdots, x_d)\to x'=(x_{m+1},x_{m+2}, \cdots, x_d)$. 
The integer $1\leq m\leq 3$ is the dimension reduction step length in our algorithm.  
In Section \ref{sec-4}, we 
shall demonstrate using numerical tests the optimal choice of step length $m$. 
 
We also recall that the tensor product quadrature rule for $I(g,\Omega)$ is defined as
\begin{equation}\label{eq2.14}
J(g, \Omega)= \sum_{i_1, i_2, \cdots, i_d=1}^N w_{i_1}w_{i_2}\cdots w_{i_d} g(\xi_{i_1,\cdots,x_{i_d}}).
\end{equation}
Where $\{\xi_j\}_{j=1}^N$ and $\{w_j\}_{j=1}^N$ are the nodes and weights of the given 1-d quadrature 
rule \eqref{eq2.6b}, and $\xi_{i_1,\cdots,x_{i_d}}=(\xi_{i_1},\xi_{i_1},\cdots, \xi_{i_d})$. Clearly, the 
computational complexity of the above formula is $O(N^d)$. 
 
Rewrite \eqref{eq2.14} as 
\begin{align}\label{eq2.15}
J(g, \Omega) &= \sum_{i_{m+1}, \cdots, i_d=1}^N w_{i_{m+1}}w_{i_2}\cdots w_{i_d}  \biggl(  \sum_{i_1, \cdots, i_m=1}^N  w_{i_1}w_{i_2}\cdots w_{i_m} g(\xi_{i_1,\cdots,x_{i_d} })  \biggr) \\
&= \sum_{i_{m+1}, \cdots, i_d=1}^N w_{i_{m+1}}w_{i_2}\cdots w_{i_d}\,  g_{d-m} \bigl(\xi_{i_{m+1}}, \cdots, \xi_{i_d} \bigr),
\nonumber
\end{align}
where 
\begin{equation}\label{eq2.14b}
g_{d-m}(s_1,\cdots,s_{d-m}) =  \sum_{i_1, \cdots, i_m=1}^N  w_{i_1}w_{i_2}\cdots w_{i_m}\, g\bigl((\xi_1, \cdots, \xi_m, s_1, \cdots, s_{d-m}) \bigr). 
\end{equation}
We note that in our MDI algorithm $g_{d-m}$ is formed as a symbolic function using \eqref{eq2.14b} and the right-hand side of \eqref{eq2.15} is a $(d-m)$-order multi-summation, which itself can be evaluated by employing the above dimension 
reduction strategy. The reduction can be iterated  $\ell:=[\frac{d}m]$ times until $d-\ell m \leq m$. 
Since $m\leq 3$, the final sum can be evaluated by calling
\textbf{Algorithm} 2.1 or 2.2. To realize this procedure, we introduce the following conventions. 

\smallskip
\begin{itemize}
	\item If $k = 1$, set MDI$(k,g_k,\Omega_k,N,m,r) :=J(g_k,\Omega_k)$, which is computed by using the one-dimensional quadrature rule \eqref{eq2.6b}. 
	\item If $k = 2$, set MDI$(k,g_k,\Omega_k,N,m,r):=$ 2d-MDI$(g_{k}, \Omega_{k}, N, r)$.
    \item If $k = 3$, set MDI$(k,g_k,\Omega_k,N,m,r):=$ 3d-MDI$(g_{k}, \Omega_{k}, N, r)$.
\end{itemize}
We note that when $k=1,2,3$, the parameter $m$ becomes a dummy variable and can be given any value. 
 
Let $P^{k-m}_{k}$ denote the natural embedding from $\mathbb{R}^k$
to $\mathbb{R}^{k-m}$ by deleting the first $m$ components of 
vectors in $\mathbb{R}^k$. Then the  tensor product quadrature approximation $J(g, \Omega)$ with $\Omega=[0,1]^d$ can be computed efficiently as follows. 
  

\begin{algorithm}[H]
	\caption{MDI($d$, $g$, $\Omega, N,m, r$) }
	\label{alg:3}
	\hspace*{0.02in} {\bf Inputs:} 
	$d (\geq 4), g, \Omega, N, m (=1,2,3),r.$ \\
	\hspace*{0.02in} {\bf Output:}
	$J=J(g,\Omega)$.
	\begin{algorithmic}[1]
		\State $\Omega_d=\Omega$, $g_d=g$, $\ell=[\frac{d}m]$.
		\For{{\color{black} $k=d:-m:d-\ell m$ (the index is decreased by $m$ at each iteration)} } 	 
		\State $\Omega_{d-m}=P_k^{k-m} \Omega_{k}$.
		\State  Construct symbolic function $g_{k-m}$ by \eqref{eq2.14c} below).
		\State  MDI$(k,g_k,\Omega_k,N,m,r)$ :=MDI$(k-m,g_{k-m},\Omega_{k-m}, N,m,r)$.
		\EndFor 
		\State $J=$ MDI$(d-\ell m, g_{d-\ell m} ,\Omega_{d-\ell m},N,m,r)$.\\
		\Return $J$.
	\end{algorithmic}
\end{algorithm}
\noindent
Where 
\begin{equation}\label{eq2.14c}
g_{k-m}(s_1,\cdots,s_{k-m}) =  \sum_{i_1, \cdots, i_m=1}^N  w_{i_1}w_{i_2}\cdots w_{i_m}\, g_k\bigl((\xi_1, \cdots, \xi_m, s_1, \cdots, s_{k-m}) \bigr). 
\end{equation}

 \begin{remark}
 	\textbf{Algorithm 2.3} recursively generates a sequence of symbolic functions $\{ g_d, g_{d-m}, g_{d-2m}, \cdots g_{d-\ell m}  \}$, each function has $m$ fewer arguments than its predecessor. As already mentioned above, our MDI algorithm 
 	explores the lattice structure of the tensor product integration points, instead of evaluating function values at all integration points independently, the MDI evaluates them in cluster and iteratively along $m$-coordinate directions, the function evaluation at any 
 	integration point is not completed until the last step of the algorithm is executed. So many computations are reused in each iteration, which is the main reason for the computation saving and to achieve a faster algorithm. Clearly, this idea can be applied to other quadrature 
 	rules, including sparse grid methods, whose integration points have a lattice-like structure.
 \end{remark}

\section{Numerical performance tests} \label{sec-3}
In this section, we shall present extensive and purposely designed numerical experiments to gauge the performance of the proposed MDI algorithm and to compare it with the standard tensor product (STP) method and the classical Monte Carlo (MC) method
for computing high dimensional integrals. All the numerical tests show that 
the MDI outperforms both TP and MC methods in low and medium dimensions (i.e., $d\approx 100$), and significantly outperforms them in very high dimensions (i.e., $d\approx 1000$), and succeeds even when the other two methods fail.  
We shall evaluate the influence of the choice of the 1-d base quadrature rule (indicated by 
the parameter $r$) and step length of the dimension iteration (indicated by 
the parameter $m$). 
 
All our numerical experiments are done in Matlab 9.4.0.813654(R2018a) on a desktop PC with Intel(R) Xeon(R) Gold 6226R CPU 2.90GHz and 32GB RAM.

\subsection{Two and three-dimensional tests} 
We first test our MDI on simple 2- and 3-d examples and to compare its performance 
(in terms of the CPU time) with the STP and MC methods. A word of warning is that  
due to small size of the problems and good accuracy of all the methods, the performance differences between of these methods 
may not be significant when the integrand $g$ is very ``nice".  This is the reason that 
we shall use oscillatory or rapidly growing integrands which often require to use a large number of integration points to achieve high accuracy.  

\medskip
{\bf Test 1.} Let $\Omega=[0,2]^2$ and consider the following 2-d integrands:
\begin{equation}\label{ex1}
 g(x):= \exp\bigl(5x_1^2+5x_2^2\bigr); \qquad \widehat{g}(x): = \sin\bigl(2\pi+10x_1^2+5x_2^2\bigr).
\end{equation}

Let $h\in (0,1)$ denote the grid size of the tensor product grid. Then the number 
of integration points in each coordinate direction is  {\color{black}$N=\frac{2}{h}+1$.} The base
1-d quadrature rule is chosen to be the Simpson's rule, hence, $r=2$. Its
composite quadrature rule is denoted by STP-S which stands for 
the standard tensor product-Simpson rule.  

Table \ref{tab:1} and \ref{tab:2} present the computational results (errors and CPU times) of the STP-S and MDI methods for approximating $I(g,\Omega)$ and $I(\widehat{g}, \Omega)$,
respectively. 
\begin{table}[H]
	\centering
	\begin{tabular}{cccccc}
		\cline{1-6} \noalign{\smallskip}
		\multicolumn{2}{c}{}&\multicolumn{2}{c}{STP-S} &\multicolumn{2}{c}{MDI} \\
		\cline{3-6} \noalign{\smallskip}
		\makecell[c]{Mesh size\\ ($h$)}&\makecell[c]{Total \\ nodes}&\makecell[c]{Relative\\ error} &\makecell[c]{CPU\\time} &\makecell[c]{Relative\\ error} &\makecell[c]{CPU\\time} \\
		\noalign{\smallskip}\hline\noalign{\smallskip}
		0.1  & 441 & $1.2146\times 10^{-1}$   & 0.0380032   &$1.2146\times 10^{-1}$ & 0.1371068 \\
		0.05  & 1681 & $1.0222\times 10^{-2}$  & 0.0438104   &$1.0222\times 10^{-2}$ &0.1857806 \\
		0.025   &6561  & $7.0238\times 10^{-4}$ & 0.0545541   &$7.0238\times 10^{-4}$    &  0.3617802  \\
		0.0125   & 25921 & $4.5031\times 10^{-5}$& 0.0633071&$4.5031\times 10^{-5}$    & 0.5514163  \\
		0.0100  &40401 & $1.8502\times 10^{-5}$ &  0.0659092 &$1.8502\times 10^{-5}$ & 0.6151638\\
		
		0.00625  &103041 &  $2.8328\times 10^{-6}$ & 0.0720637   & $2.8328\times 10^{-6}$ &0.8968891 \\
		\noalign{\smallskip}\hline
	\end{tabular}
\caption{Relative errors and CPU times of STP-S and MDI simulations with  $m=1$ for approximating $I(g,\Omega)$.}
\label{tab:1}       
\end{table}

\begin{table}[H]
	\centering
	\begin{tabular}{cccccc}
		\cline{1-6} \noalign{\smallskip}
		\multicolumn{2}{c}{}&\multicolumn{2}{c}{STP-S} &\multicolumn{2}{c}{MDI} \\
		\cline{3-6} \noalign{\smallskip}
		\makecell[c]{Mesh size\\ ($h$)}&\makecell[c]{Total \\ nodes}&\makecell[c]{Relative\\ error} &\makecell[c]{CPU\\time} &\makecell[c]{Relative\\ error} &\makecell[c]{CPU\\time} \\
		\noalign{\smallskip}\hline\noalign{\smallskip}
		0.1  & 441 & $8.4038\times 10^{-1}$   &  0.0413903   &$8.4038\times 10^{-1}$ & 0.1381629 \\
		0.05  & 1681 & $1.2825\times 10^{-2}$  & 0.0477657   &$1.2825\times 10^{-2}$ & 0.1841843 \\
		0.025   &6561  & $5.1928\times 10^{-4}$ & 0.0579602   &$5.1928\times 10^{-4}$   &  0.2845160  \\
		0.0125   & 25921 & $2.9642\times 10^{-5}$& 0.0579613&$2.9642\times 10^{-5}$   & 0.4957854  \\
		0.0100  &40401 & $1.2014\times 10^{-5}$ &  0.0617214 &$1.2014\times 10^{-5}$ & 0.6276175\\
		
		0.00625  &103041 &  $1.8123\times 10^{-6}$ & 0.0708075   & $1.8123\times 10^{-6}$ &0.9682539 \\
		0.003125  & 410881&  $1.1213\times 10^{-7}$ & 0.0946567   & $1.1213\times 10^{-7}$ &2.1355674 \\
		\noalign{\smallskip}\hline
	\end{tabular}
	\caption{Relative errors and CPU times of STP-S and MDI simulations with  $m=1$ for approximating $I(\widehat{g},\Omega)$. }	\label{tab:2} 
\end{table}

From Table \ref{tab:1} and \ref{tab:2}, we observe that the CPU times  used by these two methods are very small although that of the STP-S method in both simulations are slightly less. However,  we like to note that both methods are very efficient and the difference is almost negligible in the 2-d case. 

\medskip
{\bf Test 2.}  Let $\Omega=[0,2]^3$ and we consider the following 3-d integrands:
\begin{equation}\label{ex3}
g(x)= \exp\bigl(5x_1^2+5x_2^2+5x_3^2\bigr),\qquad \widehat{g}(x)=\sin\bigl(2\pi+10x_1^2+5x_2^2+20x_3^2\bigr).
\end{equation}
We compute integrals of these two functions over $\Omega$ by using the STP-S and MDI methods.  Again, let $h$ denote the grid size,  {\color{black}$N=\frac{2}{h}+1$,} $r=2$ and $m=1$. 

\begin{table}[H]
	\centering
	\begin{tabular}{cccccc}
	\cline{1-6} \noalign{\smallskip}
	\multicolumn{2}{c}{}&\multicolumn{2}{c}{STP-S} &\multicolumn{2}{c}{MDI} \\
	\cline{3-6} \noalign{\smallskip}
	\makecell[c]{Mesh size\\ ($h$)}&\makecell[c]{Total \\ nodes}&\makecell[c]{Relative\\ error} &\makecell[c]{CPU\\time(s)} &\makecell[c]{Relative\\ error} &\makecell[c]{CPU\\time(s)} \\
	\noalign{\smallskip}\hline\noalign{\smallskip}
	0.1  & 9261 & $1.8762\times 10^{-1}$   & 0.0594529   &$1.8762\times 10^{-1}$ & 0.1678806\\
	0.05  & 68921 & $1.0222\times 10^{-2}$   & 0.0830445   &$1.0222\times 10^{-2}$ &0.2330631 \\
	0.025   &531441  & $7.0238\times 10^{-4}$ & 0.1748109   &$7.0238\times 10^{-4}$    &  0.4138331  \\
	0.0125   & 4173281 & $4.5031\times 10^{-5}$& 0.4316260&$4.5031\times 10^{-5}$    & 0.8359030  \\
	0.0100  &8120601 & $1.8502\times 10^{-5}$ &  0.7565165 &$1.8502\times 10^{-5}$ & 1.0155151\\	
	0.00625  &33076161 &  $2.8328\times 10^{-6}$ & 2.7365100   & $2.8328\times 10^{-6}$ &1.8655724 \\
    0.003125  &263374721 &  $2.6601\times 10^{-7}$ & 56.872493   & $ 2.6601\times 10^{-7}$ &7.6742217 \\
	\noalign{\smallskip}\hline
\end{tabular}
\caption{Relative errors and CPU times of STP-S and MDI simulations  with $m=1$ for approximating $I(g,\Omega)$. }	\label{tab:3} 
\end{table}
 
\begin{table}[H]
	\centering
	\begin{tabular}{cccccc}
		\cline{1-6} \noalign{\smallskip}
		\multicolumn{2}{c}{}&\multicolumn{2}{c}{STP-S } &\multicolumn{2}{c}{MDI  } \\
		\cline{3-6} \noalign{\smallskip}
		\makecell[c]{Mesh size\\ ($h$)}&\makecell[c]{Total \\ nodes}&\makecell[c]{Relative\\ error} &\makecell[c]{CPU\\time(s)} &\makecell[c]{Relative\\ error} &\makecell[c]{CPU\\time(s)} \\
		\noalign{\smallskip}\hline\noalign{\smallskip}
		0.1  & 9261 & $2.5789\times 10^{-1}$   & 0.042875   &$2.5789\times 10^{-1}$ & 0.173581\\
		0.05  & 68921 & $3.0493\times 10^{-1}$   & 0.0729843   &$3.0493\times 10^{-1}$ &0.247261 \\
		0.025   &531441  & $1.2800\times 10^{-2}$ &  0.1942965   &$1.2800\times 10^{-2}$   &  0.4377125  \\
		0.0125   & 4173281 & $4.9563\times 10^{-4}$& 0.9138938&$4.9563\times 10^{-4}$    &  0.8930069  \\
		0.0100  &8120601 & $1.9345\times 10^{-4}$ &  1.7594922 &$1.9345\times 10^{-4}$ & 1.1012937\\	
		0.00625  &33076161 &  $2.8065\times 10^{-5}$ & 7.0588477   & $2.8065\times 10^{-5}$ &2.3235149 \\
		0.003125  &263374721 &  $1.7128\times 10^{-6}$ & 56.322503   & $1.7128\times 10^{-6}$ &9.7523139 \\
		\noalign{\smallskip}\hline
	\end{tabular}
\caption{Relative errors and CPU times of STP-S and MDI simulations with $m=1$ for approximating $I(\widehat{g},\Omega)$. }	\label{tab:4} 
\end{table}

Table \ref{tab:3} and \ref{tab:4} display the computational results (errors and CPU times), we  observe that when the number of integration points is small (i.e., the grid size $h$ is relatively large), the STP-S method requires less CPU times in both simulations. However, when the number of integration points increases, the advantage shifts to the MDI method and becomes significant when the number of integration points become large.  This is because, by the computational complexity analysis to be given in the next section,  the CPU time required by the MDI method grows in $(d^3N^2)$ order while that of the STP-S method increases in exponential order $O(N^d)$.  As a result, it is expected that when $d=3$ and  $N$ is large,
the advantage of the MDI method over any standard tensor product method becomes significant, and it will be even more pronouncing when both $d$ and $N$ become large as shown by the tests to be given in the next subsection. 

\subsection{High-dimensional tests}
Since the MDI method is designed to computing high dimensional integration, 
it is important to test its performance and power for $d>>1$. In addition, 
we provide a performance comparison 
(in terms of the CPU time) of the MDI with standard tensor product methods as well as with the classical Monte Carlo (MC) method on computing high-dimensional integration.

The next test compares the performance of the MDI and STP-S methods on computing a well-known
integral  in dimensions $2\leq d\leq 11$, respectively.

\medskip
{\bf Test 3.}  Let $\Omega=[0,1]^d$ for $2\leq d\leq 11$ and consider the following Gaussian integrand:
\begin{equation}\label{ex5}
 g(x)= \frac{1}{\sqrt{2\pi}}\exp\Bigl(-\frac{1}{2}|x|^2\Bigr),
\end{equation}
where $|x|$ stands for the Euclidean norm of the vector $x\in \mathbb{R}^d$. 

Once again, we approximate the integral $I(g,\Omega)$ by the MDI and STP-S (i.e., $r=2$) methods as done in {\bf Test 1-2}. We also set $m=1$ in the MDI method and perform the simulations with two grid sizes $h=0.1, 0.05$ (or $N=11, 21$) respectively.   

\begin{table}[H]
	\centering
	\begin{tabular}{ccccc}
		\cline{1-5} \noalign{\smallskip}
		\multicolumn{1}{c}{}&\multicolumn{2}{c}{STP-S} &\multicolumn{2}{c}{MDI} \\
		\cline{2-5} \noalign{\smallskip}
		\makecell[c]{Dimension\\($d$)}&\makecell[c]{Relative\\ error} &\makecell[c]{CPU\\time(s)} &\makecell[c]{Relative\\ error} &\makecell[c]{CPU\\time(s)} \\
		\noalign{\smallskip}\hline\noalign{\smallskip}
		2   & $1.5809\times 10^{-6}$ &  0.0015542 & $1.5809\times 10^{-6}$  &0.0853444  \\
		4   & $3.1618\times 10^{-6}$ & 0.0091310 &$3.1618\times 10^{-6}$    &0.1348654  \\
		6 & $4.7427\times 10^{-6}$ & 0.4403814 & $4.7427\times 10^{-6}$   & 0.5389767   \\
		8  & $6.3237\times 10^{-6}$ & 56.1856842  &$6.3237\times 10^{-6}$    &1.4880431   \\
		10 & $7.9046\times 10^{-6}$  &  7341.3815698&$7.9046\times 10^{-6}$    & 3.7304532  \\  
		11 & $8.6951\times 10^{-6}$ &80322.5805531   & $8.6951\times 10^{-6}$    &  5.2628807  \\  
		\noalign{\smallskip}\hline
	\end{tabular}
\caption{Relative errors and CPU times of STP-S and MDI simulations with $m=1$ for approximating $I(g,\Omega)$ when $N=11$. }	\label{tab:55}     
\end{table}
\begin{table}[H]
	\centering
	\begin{tabular}{ccccc}
		\cline{1-5} \noalign{\smallskip}
		\multicolumn{1}{c}{}&\multicolumn{2}{c}{STP-S} &\multicolumn{2}{c}{MDI} \\
		\cline{2-5} \noalign{\smallskip}
		\makecell[c]{Dimension\\($d$)}&\makecell[c]{Relative\\ error} &\makecell[c]{CPU\\time(s)} &\makecell[c]{Relative\\ error} &\makecell[c]{CPU\\time(s)} \\
		\noalign{\smallskip}\hline\noalign{\smallskip}
		2   & $9.8542\times 10^{-8}$ &  0.0068182 & $9.8542\times 10^{-8}$  &0.1122873  \\
		4   & $1.9708\times 10^{-7}$ & 0.0715035 &$1.9708\times 10^{-7}$    &0.6049681  \\
		6 & $2.9564\times 10^{-7}$ & 20.7002115 & $2.9564\times 10^{-7}$   & 6.4742624   \\
		8  & $3.9149\times 10^{-7}$ & 9622.1118103  &$3.9149\times 10^{-7}$    &19.8850829    \\
		9 & $4.4344\times 10^{-7}$ &  215136.0654597   & $4.4344\times 10^{-7}$   &28.6823906   \\ 
		10 & $4.9271\times 10^{-7}$  & failed &$4.9271\times 10^{-7}$    & 38.9745044  \\   
		\noalign{\smallskip}\hline
	\end{tabular}
	\caption{Relative errors and CPU times of STP-S and MDI simulations with  $m=1$ for approximating $I(g,\Omega)$ when $N=21$.}\label{tab:5}  
\end{table} 
\begin{figure}[H]
	\centerline{
	 \includegraphics[width=2.5in, height=1.8in]{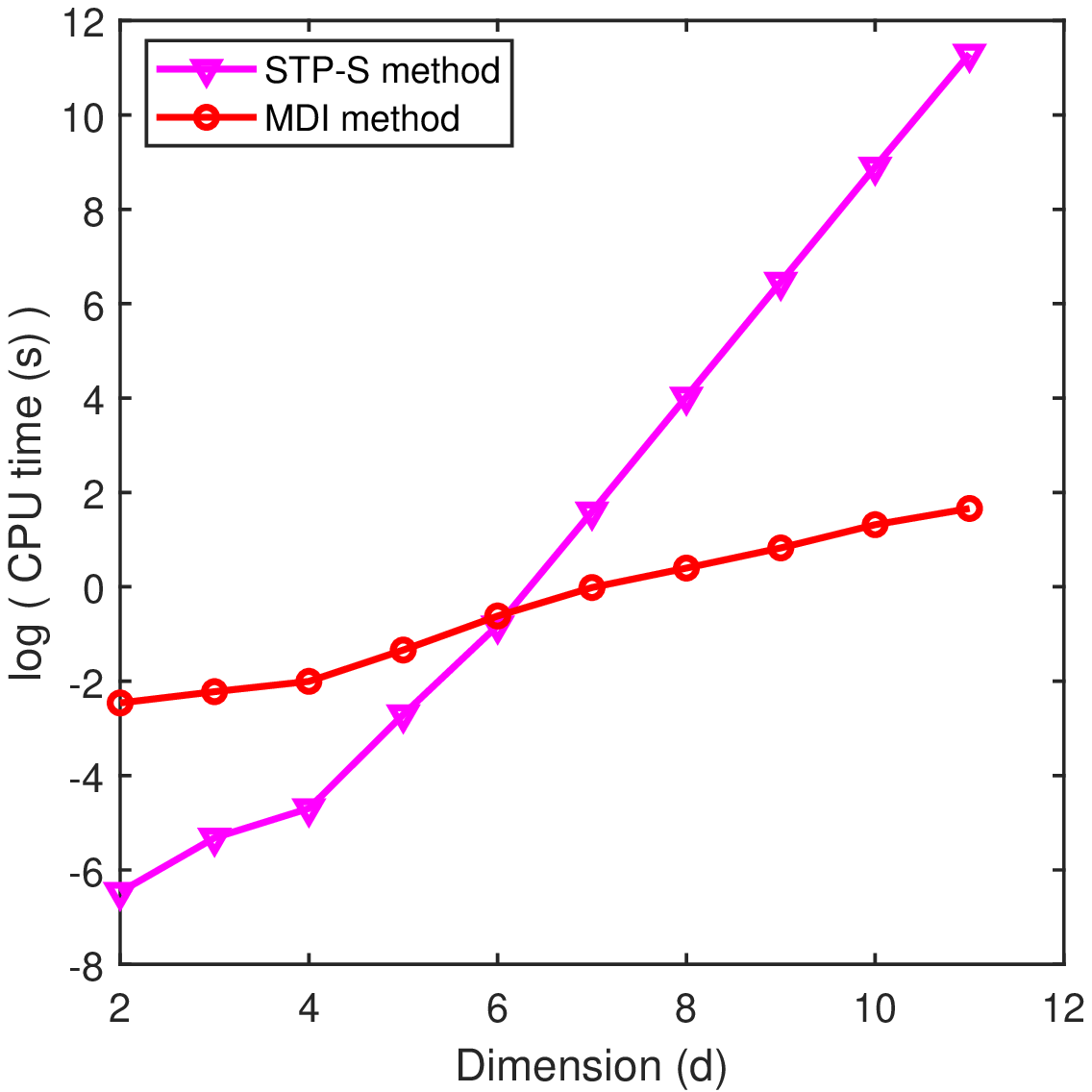}
	 \includegraphics[width=2.5in, height=1.8in]{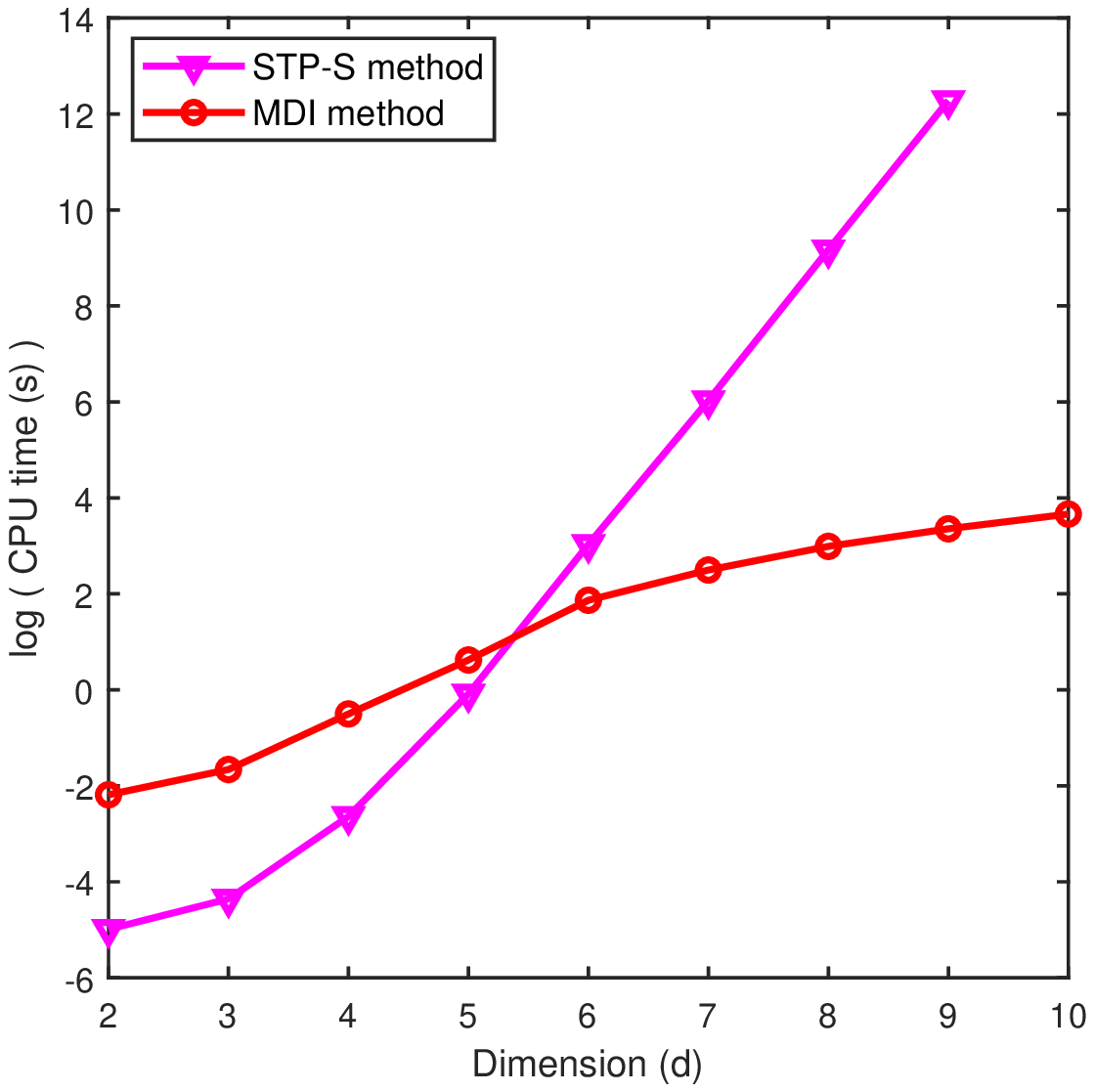}
	}
	\caption{CUP time comparison of STP-S and MDI simulations: $N=10$ (left), $N=20$ (right).}
	\label{fig.1}       
\end{figure}

Table \ref{tab:55} presents the relative errors and CPU times of both MDI and STP-S
methods using the grid size $h=0.1$ (or $N=11$) and Table \ref{tab:5} gives the corresponding results for $h=0.05$ (or $N=21$).
We observe from both tables that the errors of both methods are the same as they should be (because they compute the same multi-summation in each simulation), but their CPU times are significantly  different. The STP-S method is more efficient when the dimension $d\leq 6$ when $N=11$ and $d\leq 5$ when $N=21$, but the MDI method excels when the dimension $d>6$ and the winning margin becomes  significant as $d$ and $N$ increase (also see Figure \ref{fig.1}).  For example, when $d=11$ and $N=11$, the CPU time required by the STP-S method is about $80323$ seconds, which is about $22$ hours, but the CPU time required by the MDI method is only less than $6$ seconds!  In addition, 
when $d=10$ and $N=21$, the STP-S method fails to compute the integral (because the 
computational cost is too large for the computer to handle), but it only takes the MDI method about $39$ seconds to finish the computation! 
The reason for such a dramatic CPU time saving is,  by the computational complexity analysis to be given in the next section,  that the CPU time required by the MDI method grows in $(d^3N^2)$ order while that of the STP-S method increases in exponential order $O(N^d)$. 


\smallskip
Arguably, the classical Monte Carlo (MC) method (or its variant) is the method of choice for computing high dimensional integration. However, due to its low order of convergence and intrinsic need for using large amount of samples,  it is only capable of simulating low and medium dimensional integration 
in practice due to the large number of function evaluations at randomly sampled integration points, which also grows quickly as the dimension $d$ increases {\color{black} (due to the rapid growth of the variance)}.  In the next test, we compare the performance of the MDI (with parameters $r=2,N=11,m=1$) and the classical MC method.

\medskip
{\bf Test 4.} Let $\Omega=[0,1]^d$ and choose the following integrands: 
\begin{equation}\label{ex6}
	 g(x)= \prod_{i=0}^{d} \frac{1}{0.9^2+(x_i-0.6)^2}, 
	 \qquad  \widehat{g}(x)= \frac{1}{\sqrt{2\pi}}\exp\Bigl(-\frac{1}{2}|x|^2\Bigr).
\end{equation}

First, we use the relative error as the metric to guage the performance, namely, we use sufficient number of the random sampling points for the MC method so it produces a comparable relative error to that of the MDI method. The computational results for   
approximating $I(g,\Omega)$ and $I(\widehat{g}, \Omega)$ are 
presented respectively in Table \ref{tab:6} and \ref{tab:7} below. 
\begin{table}[H]
	\centering
	\begin{tabular}{ccccc}
		\cline{1-5} \noalign{\smallskip}
\multicolumn{1}{c}{}&\multicolumn{2}{c}{MC} &\multicolumn{2}{c}{MDI} \\
\cline{2-5} \noalign{\smallskip}
\makecell[c]{Dimension\\($d$)}&\makecell[c]{Relative\\ error} &\makecell[c]{CPU\\time(s)} &\makecell[c]{Relative\\ error} &\makecell[c]{CPU\\time(s)} \\
\noalign{\smallskip}\hline\noalign{\smallskip}

        5  & $2.6251\times 10^{-5}$ & 9.2872030  & $2.6251\times 10^{-5}$   & 0.1755604  \\
		10&$5.2504\times 10^{-5}$   &17.9837795   & $5.2504\times 10^{-5}$ &0.2510754  \\
		20 &$1.0501\times 10^{-4}$   &66.6797401    & $1.0501\times 10^{-4}$ &0.6229516   \\ 
        30  &$ 1.5752\times 10^{-4}$   &4302.1801130  & $ 1.5752\times 10^{-4}$ & 0.8786786    \\  
        35  &$1.8377\times 10^{-4}$   &11055.6351555  & $1.8377\times 10^{-4}$ & 1.0470143    \\ 
		40 &$2.1003\times 10^{-4}$  &  failed & $2.1003\times 10^{-4}$ & 1.2647430   \\ 
		80 &$4.2011\times 10^{-4}$   &  & $4.2011\times 10^{-4}$ &2.8572520     \\ 
		100 & $5.2516\times 10^{-4}$ &  & $5.2516\times 10^{-4}$ &3.2840268   \\ 
		\noalign{\smallskip}\hline
	\end{tabular}
	\caption{CPU times of the MDI and MC  simulations with  comparable relative errors for approximating $I(g, \Omega)$.} \label{tab:6}     
\end{table}
\begin{table}[H]
	\centering
	\begin{tabular}{ccccc}
		\cline{1-5} \noalign{\smallskip}
\multicolumn{1}{c}{}&\multicolumn{2}{c}{MC} &\multicolumn{2}{c}{MDI} \\
\cline{2-5} \noalign{\smallskip}
\makecell[c]{Dimension\\($d$)}&\makecell[c]{Relative\\ error} &\makecell[c]{CPU\\time(s)} &\makecell[c]{Relative\\ error} &\makecell[c]{CPU\\time(s)} \\
\noalign{\smallskip}\hline\noalign{\smallskip}
5  & $3.9523\times 10^{-6}$ &33.3246187  & $3.9523\times 10^{-6}$ & 0.2620139  \\
10&$7.9046\times 10^{-6}$   &346.0491007   &$7.9046\times 10^{-6}$    & 3.7304532  \\ 
20 &$1.5809\times 10^{-5}$   &1210.092329    & $1.5809\times 10^{-5}$ &35.9158247    \\ 
30  &$2.3714\times 10^{-5}$   &3416.6898735  & $2.3714\times 10^{-5}$ &108.7032555    \\
35 &$2.7666\times 10^{-5}$   & 12664.1000000   & $2.7666\times 10^{-5}$ &154.1259392    \\   
40 &$3.1618\times 10^{-5}$   & failed & $ 3.1618\times 10^{-5}$ &240.1254540   \\ 
80 &$6.3238\times 10^{-5}$   &  & $6.3238\times 10^{-5}$ &1678.4531292     \\ 
100 &$7.9049\times 10^{-5}$   &  & $7.9049\times 10^{-5}$ &3207.3510623   \\ 
		\noalign{\smallskip}\hline
	\end{tabular}
\caption{CPU times of the MDI and MC simulations with comparable relative errors for approximating $I(\widehat{g}, \Omega)$.} \label{tab:7}  
\end{table}

From Table \ref{tab:6} and \ref{tab:7}, we clearly see that the CPU times of the MDI  and MC methods are significantly different, the discrepancy becomes so dramatic when $d\geq 40$ because the MDI method only takes a few seconds/minutes to finish the computation
of approximating $I(g,\Omega)$ and $I(\widehat{g},\Omega)$, while the MC method fails to produce a result on the computer! This is because, in order for the MC method 
to obtain an approximate value with the relative error of order $10^{-5}$, it requires  about $10^{10}$ randomly sampled integration points at which function values must be computed independently, a task that is too big to be handled by Matlab on the computer. 
We note that $g$ is an oscillatory function and $\widehat{g}$ is an exponentially 
growth function, both functions are tough for the MC method to handle,
because a very large number of sampling points must be used to resolve those functions with a reasonable resolution and the function values must be computed independently at those points in the MC method. On the other hand, 
although the MDI method must use a comparable large amount (if not larger) of integration points (because the underlying tensor product method does), 
due to its efficient way of computing those function evaluations in cluster and iteratively along each coordinate direction,
 the MDI can handle the computation of the multi-summation in the blink of an eye. 

\begin{table}[htb]
	\centering
	\begin{tabular}{cccccc}
		\cline{1-6} \noalign{\smallskip}
		\multicolumn{1}{c}{}&\multicolumn{1}{c}{}&\multicolumn{2}{c}{MC} &\multicolumn{2}{c}{MDI} \\
		\cline{3-6} \noalign{\smallskip}
		\makecell[c]{Dimension\\($d$)}&\makecell[c]{Total\\nodes}&\makecell[c]{Relative\\ error} &\makecell[c]{CPU\\time(s)} &\makecell[c]{Relative\\ error} &\makecell[c]{CPU\\time(s)} \\
		\noalign{\smallskip}\hline\noalign{\smallskip}
5  &$11^5$& $5.8241\times 10^{-4}$ &0.3056289  & $2.6251\times 10^{-5}$   & 0.1755604  \\
11&$11^{11}$&$4.3309\times 10^{-7}$  &9827.3992235   & $5.7754\times 10^{-5}$ &0.2815804 \\
15&$11^{15}$& failed  &  & $7.8757\times 10^{-5}$ &  0.4811588 \\
30  &$11^{30}$&   &  & $ 1.5752\times 10^{-4}$ & 0.8786786    \\  
40 &$11^{40}$&  &   & $2.1003\times 10^{-4}$ & 1.2647430   \\ 
50 &$11^{50}$&  &   & $2.6254\times 10^{-4}$ &1.6818385  \\ 
\noalign{\smallskip}\hline
	\end{tabular}
	\caption{CPU times of the MDI and MC simulations using the same number of integration points for approximating $I(g,\Omega)$.}
\label{tab:66}
\end{table}

Next, we compute the same test problems as above but use a different metric to gauge 
the performance of both methods. We now specify the number of integration
points, instead of the relative errors, then let both methods compute their 
respective multi-summations using the same number of points (so the same number of function evaluations are required by both methods to approximate the integrals). 
Table \ref{tab:66} and \ref{tab:77} present the simulation results. Although the 
numbers are slightly different but the message is the same, that is, the MC method 
fails to produce a result when the dimension $d\geq 15$ while the MDI can finish
the simulation in the blink of an eye for the integrand $g$ and in a few minutes 
for the integrand $\widehat{g}$ even when $d=50$.

\begin{table}[htb]
	\centering
	\begin{tabular}{cccccc}
		\cline{1-6} \noalign{\smallskip}
		\multicolumn{1}{c}{}&\multicolumn{1}{c}{}&\multicolumn{2}{c}{MC} &\multicolumn{2}{c}{MDI} \\
		\cline{3-6} \noalign{\smallskip}
		\makecell[c]{Dimension\\($d$)}&\makecell[c]{Total\\ nodes}&\makecell[c]{Relative\\ error} &\makecell[c]{CPU\\time(s)} &\makecell[c]{Relative\\ error} &\makecell[c]{CPU\\time(s)} \\
		\noalign{\smallskip}\hline\noalign{\smallskip}
5  &$11^{5}$&$1.0664\times 10^{-3}$  &0.0080803  & $3.9523\times 10^{-6}$ & 0.2620139  \\
11&$11^{11}$&$1.1947\times 10^{-6}$  &10195.085484 &$8.6951\times 10^{-6}$    & 6.9215447  \\ 
15&$11^{15}$& failed  &   &$1.1856\times 10^{-5}$    & 15.316679  \\ 
30  &$11^{30}$&   &  & $2.3714\times 10^{-5}$ &108.703255    \\
40 &$11^{40}$&   &    & $ 3.1618\times 10^{-5}$ &240.125454   \\ 
50 &$11^{50}$&   &   & $3.9523\times 10^{-5}$ &413.607179    \\ 
		\noalign{\smallskip}\hline
	\end{tabular}
	\caption{CPU times of the MDI and MC simulations using the same number of integration points for approximating $I(\widehat{g},\Omega)$.}
\label{tab:77}
\end{table}

\medskip
One natural question is how high the dimension $d$ which the MDI can handle. 
First, we note that the answer is machine-dependent as expected. Second, we 
perform the next test to seek an answer to this question using the computer 
at our disposal as described at the beginning of this section. 

\medskip
{\bf Test 5.} Let $\Omega=[0,1]^d$ and consider the following integrands:

\begin{equation}\label{ex-5}
g(x)= \exp\Bigl(\sum_{i=1}^{d}(-1)^{i+1}x_i \Bigr), 
\qquad \widehat{g}(x)= \prod_{i=0}^{d} \frac{1}{0.9^2+(x_i-0.6)^2}.
\end{equation}
We then  approximate $I(g,\Omega)$ and $I(\widehat{g},\Omega)$ using {\color{black} 
the MDI algorithm  with parameters  {\color{black}$r=2$ (composite Simpson's rule),  $m=1$, $N=7$}} and an increasing sequence of $d$. 
The computational 
results are presented in Table \ref{tab:17} and \ref{tab:16}. The simulation is stopped  at $d=1000$ because it is already in the very high
dimension regime and $N=7$ is chosen to minimize the computation and because
it is sufficient to produce reasonable relative errors. This test demonstrates the 
promise and capability of the MDI method for efficiently computing high dimensional integrals.

\begin{table}[H]
	\centering
	\begin{tabular}{cccc}
		\toprule
		\makecell[c]{Dimension\\($d$)} &\makecell[c] {Total \\ nodes }  &\makecell[c]{Relative\\error}&\makecell[c]{CPU\\time(s)}\\
		
		\midrule
		10& $7^{10}$ & $4.2726\times 10^{-5}$ & 0.2459398     \\  
		100&$7^{100}$ & $4.2734\times 10^{-4}$ & 74.6973942     \\ 
		200 &$7^{200}$& $8.5487\times 10^{-4}$  & 503.3034372    \\ 
		300&$7^{300}$ &$1.4386\times 10^{-4}$  &  1560.0488058    \\ 
		400 &$7^{400}$&$ 1.7097\times 10^{-3}$  &  3546.4398972     \\ 
		500 &$7^{500}$& $2.1371\times 10^{-3}$   &  6772.0225935   \\ 
		600 &$7^{600}$& $2.8772\times 10^{-3}$   &  11954.1886240    \\ 
		700 &$7^{700}$& $3.3566\times 10^{-3}$   & 19355.4847153    \\ 
		800 &$7^{800}$& $3.5194\times 10^{-3}$   & 28273.5752793  \\ 
		900 &$7^{900}$& $3.8467\times 10^{-3}$   &  42427.2391457     \\ 
		1000&$7^{1000}$ & $4.2742\times 10^{-3}$   &  62445.0882189    \\ 
		\noalign{\smallskip}\hline
	\end{tabular}
\caption{Computation results for approximating $I(g,\Omega)$ by the  {\color{black}MDI} algorithm.} \label{tab:17}     
\end{table}
\begin{table}[H]
	\centering
		\begin{tabular}{cccc}
			\toprule
			\makecell[c]{Dimension\\($d$)} &\makecell[c] {Total\\nodes}  &\makecell[c]{Relative\\error}&\makecell[c]{CPU\\time(s)}\\
			
			\midrule
			10    &$7^{10}$ & $4.0743\times 10^{-4}$ &0.2059168 \\  
			100   & $7^{100}$ & $4.0818 \times 10^{-3}$ & 2.0993900 \\ 
			200   & $7^{200}$   & $8.1803 \times 10^{-3}$ & 4.5213100   \\
			300   &  $7^{300}$    & $1.2295\times 10^{-2}$ & 7.2644682 \\
			400   & $7^{400}$   & $1.6427 \times 10^{-2}$ &10.1062101    \\
			500   & $7^{500}$  & $2.0576 \times 10^{-2}$ &13.5705851   \\
			600   & $7^{600}$   & $2.4742 \times 10^{-2}$ & 17.8284828   \\
			700   & $7^{700}$   & $2.8925 \times 10^{-2}$ & 21.6876065   \\
			800   & $7^{800}$   & $3.3125 \times 10^{-2}$ & 25.9204244    \\
			900   & $7^{900}$    & $3.7342\times 10^{-2}$ &31.3307727   \\
			1000   & $7^{1000}$   & $4.1576 \times 10^{-2}$ &35.7704489   \\
			\bottomrule
	\end{tabular}
	\caption{Computation results for approximating $I(\widehat{g},\Omega)$ by the  {\color{black}MDI} algorithm.}\label{tab:16}
\end{table}

\section{Influence of parameters}\label{sec-4}
Besides the dimension $d$, there are three other input parameters in the MDI algorithm, they are $r,m$ and $N$. The parameter $r$ indicates the choice of 1-d base numerical 
quadrature rule. As mentioned earlier, here we only consider four such choices, 
hence, $r$ takes integer values $\{1,2,3,4\}$ and they represent respectively the (composite) trapezoidal rule, Simpson's rule, two-point Gaussian rule, and midpoint rule. Their efficiency will be tested 
in this section. 

Recall that $m$ represents the step length in the multi-dimension iteration, namely, it indicates how many dimensions to reduce at each iteration. Practically, $1\leq m\leq 3$, hence, it takes integer values $\{1, 2, 3\}$. The performance of each of these choices will be compared in this section. 
{\color{black} It should be noted that after $\ell:=[\frac{d}m]$ iterations, the residual dimension satisfies $d-\ell m\leq m$. Then in case $m=2$ or $3$, one has two options to choose to complete the 
	algorithm. On one hand, one just continues the dimension reduction by calling 3d-MDI or 2d-MDI as explained in the definition of \textbf{Algorithm 2.3}. 
On the other hand, it is also possible to compute the remaining 2- or 3-d integration directly using the underlying 2- or 3-d tensor product formula without further dimension reduction. The effect of these two choices will be tested in this section. 

It is clear that the larger $N$, the more expensive the computation.  The dependence of the efficiency of the MDI algorithm on the parameter $N$ will also be tested.
}

\subsection{Influence of parameter $r$}
We first examine the effect of the choices $r=1,2,3,4$ in the MDI algorithm. 
They will be done on the same grid (i.e., $N$ fixed) and with the same step length $m=1$. 

\medskip 
{\bf Test 6:} Let $\Omega=[0,1]^d$ and the integrand $g$ be given by \eqref{ex5}.  
\begin{table}[h]
	\centering
	\setlength{\tabcolsep}{4mm}{
		\begin{tabular}{ccccc}
			\toprule
			\makecell[c]{Parameter\\($r$)} &\makecell[c]{Dimension\\($d$)}&\makecell[c] {Points\\($N$)}  &\makecell[c]{Relative\\error}&\makecell[c]{CPU\\time(s)}\\
			
			\midrule
			\multirowcell{10}{$r=1$}&10&11&$5.8935\times 10^{-3}$& 3.8294275 \\
			&30&11&$1.7576\times 10^{-2}$ &97.3001349\\
			&50&11&$2.9122\times 10^{-2}$ &419.1940599\\
			&70&11&$4.0532\times 10^{-2}$ &1056.501204\\
			&90&11&$5.1808\times 10^{-2}$ & 2199.34273\\
			&100&11&$5.7396\times 10^{-2}$ &3255.997642\\
			
			\midrule
			\multirowcell{10}{$r=2$}&10&11&$7.9046\times 10^{-6}$    & 3.7304532  \\
			&30&11& $2.3714\times 10^{-5}$ &108.703255    \\
			&50&11& $3.9523\times 10^{-5}$ &413.607179    \\
			&70&11& $5.5333\times 10^{-5}$ &1147.446169    \\ 
			&90&11& $7.1144\times 10^{-5}$ & 2388.382073   \\ 
			&100&11& $7.9049\times 10^{-5}$ &3207.351062   \\
		 
			\midrule
			\multirowcell{10}{$r=4$}&10&10&$2.9593\times 10^{-3}$&1.6020058\\
			&30&10&$8.9042\times 10^{-3}$ &50.2901972\\
			&50&10&$1.4884\times 10^{-2}$ & 186.0984739\\
			&70&10&$2.0899\times 10^{-2}$ &436.3740376\\
			&90&10&$2.6951\times 10^{-2}$ &855.0901709\\
			&100&10&$2.9990\times 10^{-2}$ &1062.3387568\\
			\bottomrule
	\end{tabular}}
	\caption{Efficiency comparison of the MDI algorithm with $m=1$ and  $r=1,2,4$.} \label{tab:12}
\end{table}

Table \ref{tab:12} presents the simulation results of {\bf Test 6}.
We note that since the composite two-point Gaussian rule ($r=3$) is too expensive to compute this integral when the dimension 
is larger than $10$, so it is not included in this test.  It shows that Simpson's and trapezoidal rules have the same efficiency, but Simpson's rule has 
much better accuracy. The midpoint and trapezoidal rules have the same accuracy, but the midpoint rule is 
three times more efficient than the trapezoidal rule in this test. 

\medskip 
{\bf Test 7:} Let $\Omega=[0,1]^d$ and choose the integrand $g$ as 
\begin{equation}\label{ex-9}
	g(x)= \exp\Bigl(\sum\limits_{i=1}^{d}(-1)^{i+1}x_i \Bigr), \qquad 
	\widehat{g}(x)= \prod_{i=0}^{d} \frac{1}{0.9^2+(x_i-0.6)^2} .
\end{equation}

\begin{table}[ht]
	\centering
	\setlength{\tabcolsep}{4mm}{
		\begin{tabular}{ccccc}
			\toprule
			\makecell[c]{Parameter\\($r$)} &\makecell[c]{Dimension\\($d$)}&\makecell[c] {Points\\($N$)}  &\makecell[c]{Relative\\error}&\makecell[c]{CPU\\time(s)}\\
			 
			 \midrule
			 \multirowcell{10}{$r=1$}&10&11&$8.3632\times 10^{-3}$& 0.3793491 \\
			 &30&11&$2.5300\times 10^{-2}$ &3.9242551\\
			 &50&11&$4.2521\times 10^{-2}$ & 20.7387841\\
			 &70&11&$6.0032\times 10^{-2}$ &76.7165061\\
			 &90&11&$7.7837\times 10^{-2}$ &170.2491139\\
			 &100&11&$8.6851\times 10^{-2}$ &234.4891902\\
			 
			\midrule
			\multirowcell{10}{$r=2$}&10&11&$5.5489\times 10^{-6}$    & 0.3435032  \\
			&30&11& $1.6646\times 10^{-5}$ & 4.0590394   \\
			&50&11& $2.7745\times 10^{-5}$ &20.6479181    \\
			&70&11& $3.8843\times 10^{-5}$ &69.2068795    \\ 
			&90&11& $4.9941\times 10^{-5}$ & 162.1716159   \\ 
			&100&11& $5.5491\times 10^{-5}$ & 209.2587748   \\
			 			
			\midrule
			\multirowcell{6}{$r=3$}&10&6& $2.8477\times 10^{-5}$ &0.8652789\\
			&30&6& $8.5428\times 10^{-5}$ &219.3882758  \\
			&50&6& $1.4237\times 10^{-4}$ &5281.6063020  \\
			&60&6& $2.8775\times 10^{-5}$ &16366.6127593  \\
			&70&6&   & failed  \\

			\midrule
			\multirowcell{10}{$r=4$}&10&10&$4.1576\times 10^{-3}$&0.3168370\\
			&30&10&$1.2421\times 10^{-2}$ &2.3652108\\
			&50&10&$2.0616\times 10^{-2}$ & 13.5012784\\
			&70&10&$2.8743\times 10^{-2}$ &45.0645030\\
			&90&10&$3.6802\times 10^{-2}$ &112.2156751\\
			&100&10&$4.0807\times 10^{-2}$ &153.9030288\\
			\bottomrule
	\end{tabular}}
\caption{Efficiency comparison of the MDI algorithm with $m=1$ and $r=1,2,3,4$.} \label{tab:100}
\end{table}
\begin{table}[ht]
	\centering
	\setlength{\tabcolsep}{4mm}{
		\begin{tabular}{ccccc}
			\toprule
			\makecell[c]{Parameter\\($r$)} &\makecell[c]{Dimension\\($d$)}&\makecell[c] {Points\\($N$)}  &\makecell[c]{Relative\\error}&\makecell[c]{CPU\\time(s)}\\
			
			\midrule
			\multirowcell{10}{$r=1$}&10&11&$1.2809\times 10^{-2}$&  0.2997928 \\
			&30&11&$3.7939\times 10^{-2}$ &0.8831823\\
			&50&11&$6.2429\times 10^{-2}$ &1.5778940\\
			&70&11&$8.6296\times 10^{-2}$ &2.2930658\\
			&90&11&$1.0955\times 10^{-1}$ &3.1440537 \\
			&100&11&$1.2096\times 10^{-1}$ &3.5166433\\

			\midrule
			\multirowcell{10}{$r=2$}&10&11& $5.2504\times 10^{-5}$ &0.2510754   \\
			&30&11& $ 1.5752\times 10^{-4}$ & 0.8786786       \\
			&50&11& $2.6254\times 10^{-4}$ &1.6818385   \\
			&70&11& $3.6758\times 10^{-4}$ &2.2283636   \\ 
			&90&11& $4.7263\times 10^{-4}$ & 3.1751084   \\ 
			&100&11& $5.2516\times 10^{-4}$ &3.2840268   \\
		 
			\midrule
			\multirowcell{10}{$r=3$}&10&10& $3.5014\times 10^{-5}$ &0.3039294\\
			&30&10& $1.0503\times 10^{-4}$ &1.3949243  \\
			&50&10& $1.7505\times 10^{-4}$ &4.2661412  \\
			&70&10& $2.4507\times 10^{-4}$ &9.8203358\\
			&90&10& $3.1508\times 10^{-4}$ &19.3491396\\
			&100&10& $3.5008\times 10^{-4}$ &25.3835608\\
			
			\midrule
			\multirowcell{10}{$r=4$}&10&10&$6.4658\times 10^{-3}$& 0.2454473\\
			&30&10&$1.9523\times 10^{-2}$ &0.7190225\\
			&50&10&$3.2750\times 10^{-2}$ &1.3526353 \\
			&70&10&$4.6148\times 10^{-2}$ &1.8832556\\
			&90&10&$5.9720\times 10^{-2}$ &2.4435807\\
			&100&10&$6.6572\times 10^{-2}$ &2.8005838\\
			\bottomrule
	\end{tabular}}
\caption{Efficiency comparison of the MDI algorithm with $m=1$ and $r=1,2,3,4$.} \label{tab:10}
\end{table}

Table \ref{tab:100} presents the simulation results of {\bf Test 7} for approximating integral $I(g, \Omega)$. Again, choosing different types of the 1-d base quadrature rule has a significant impact on the accuracy and 
efficiency of the MDI algorithm. 
In terms of accuracy, the trapezoidal ($r=1$) and midpoint ($r=4$) rules are comparable, but the midpoint rule is more efficient (in terms of the CPU time) as the dimension $d$ increases. 
Similarly, Simpson's ($r=2$) and two-point Gaussian ($r=3$) rules are comparable in accuracy, but Simpson's  rule is significantly more efficient even the Gaussian rule uses fewer integration points. Moreover, Simpson's 
rule is much more accurate than the trapezoidal and midpoint rules with comparable efficiency 
{\color{black} because all three quadrature rules use the same number of integration points.} The comparison shows 
that Simpson's rule is a clear winner among these four rules when they are used as the 
building block in the MDI algorithm for high dimension integration. We note that the reason 
that the two-point Gaussian rule requires a lot more CUP time is because it is costly to generate 
the Gauss points on fly and to do their function evaluations. 

Table \ref{tab:10} shows the simulation results of {\bf Test 7} for approximating integral $I(\widehat{g}, \Omega)$. Due to the nicer behavior of the integrand $\widehat{g}$, the MDI algorithm is very fast 
with all four base quadrature rules for computing this integral although the same 
observations as above can be made. Once again, the Simpson's rule excels. 

\subsection{Influence of parameter $m$}
Recall that $m$ stands for the step length in the MDI algorithm, it represents how many dimensions
are reduced at each iteration.  The intuition is the more reduction the better. However, that is not 
true because at each iteration, many $m$-dimensional tensor product sums must be evaluated. Hence, practically
we have $1\leq m\leq 3$. The next test presents a performance comparison of the MDI algorithm using 
$m=1,2,3$.

\medskip
{\bf Test 8.} Let $\Omega=[0,1]^d$, $g$ and $\widehat{g}$ be the same as in \eqref{ex6}.  

\medskip
We compute these integrals using the MDI algorithm with  {\color{black}$r=2$ (composite Simpson's rule) and $N=11$}.
Table \ref{tab:21} and \ref{tab:11} present respectively the computed results
for these two integrals. We observe that the MDI algorithm with different parameters $m$ has the same accuracy
which is expected. However, the choice of $m$ do affect the efficiency of the algorithm. It shows that the algorithm is most efficient when $m=1$. The explanation for this observation is that when using larger $m$, the number of nested   loops increases despite the number of iterations decreases. When $m=1$,  there is only one loop per iteration, so the MDI algorithm becomes faster.  

\begin{table}[H]
	\centering
	\begin{tabular}{ccccc}
		\cline{1-5} \noalign{\smallskip}
		\multicolumn{2}{c}{}&\multicolumn{1}{c}{$m=1$} &\multicolumn{1}{c}{$m=2$} &\multicolumn{1}{c}{$m=3$} \\
		\cline{3-5} \noalign{\smallskip}
		\makecell[c]{Dimension\\($d$)}&\makecell[c]{Relative\\ error} &\makecell[c]{CPU\\time(s)}  &\makecell[c]{CPU\\time(s)} &\makecell[c]{CPU\\time(s)}  \\
		\noalign{\smallskip}\hline\noalign{\smallskip}
		10   & $5.2504\times 10^{-5}$ &0.2510754     & 1.7061606    &10.8277897 \\
		30     & $ 1.5752\times 10^{-4}$ & 0.8786786    &5.6306318 &49.1959065   \\  
		50    & $2.6254\times 10^{-4}$ &1.6818385  &10.1213779  &79.0436736  \\ 
		70   & $3.6758\times 10^{-4}$ &2.2283636  &14.4966060   &127.4604677    \\ 
		90   & $4.7263\times 10^{-4}$ & 3.1751084     &19.3119927 &171.5185236 \\ 
		100   & $5.2516\times 10^{-4}$ &3.2840268   & 21.8707585 & 196.2232037  \\ 
		\noalign{\smallskip}\hline
	\end{tabular}
\caption{Efficiency comparison of the MDI algorithm with  {\color{black}$r=2$} and $m=1,2,3$.}
\label{tab:21}      
\end{table}

\begin{table}[H]
	\centering
	\begin{tabular}{ccccc}
		\cline{1-5} \noalign{\smallskip}
	\multicolumn{2}{c}{}&\multicolumn{1}{c}{$m=1$} &\multicolumn{1}{c}{$m=2$} &\multicolumn{1}{c}{$m=3$} \\
	\cline{3-5} \noalign{\smallskip}
	\makecell[c]{Dimension\\(d)}&\makecell[c]{Relative\\ error} &\makecell[c]{CPU\\time(s)}  &\makecell[c]{CPU\\time(s)} &\makecell[c]{CPU\\time(s)}  \\
	\noalign{\smallskip}\hline\noalign{\smallskip}
		10   &$7.9046\times 10^{-6}$    & 3.7304532    &17.4647858 &31.7445041 \\
		30     & $2.3714\times 10^{-5}$ &108.703255  &548.4820562 &3762.781656   \\  
		50    & $3.9523\times 10^{-5}$ &413.607179   &2208.4417154 &15367.76577  \\ 
		70   & $5.5333\times 10^{-5}$ &1147.446169  &4332.8074455 &35433.0000001    \\ 
		90   & $7.1144\times 10^{-5}$ & 2388.382073    &12397.0363557 &75097.6376617 \\ 
		100   & $7.9049\times 10^{-5}$ &3207.351062   & 16428.9713811 &102930.139707  \\ 
		\noalign{\smallskip}\hline
	\end{tabular}
\caption{Efficiency comparison of the MDI algorithm with {\color{black} $r=2$ }and $m=1,2,3$.}
\label{tab:11}    
\end{table}

\subsection{Influence of the parameter $N$} \label{sec-4.3} 
In this section, we test the influence of the number of integration points $N$ in each coordinate direction on the MDI algorithm. We set $m=1$ and $r=2$ (Simpson) in the test. 

\medskip 
{\bf Test 9.} Let $\Omega=[0,1]^d$ and choose the following integrands:
 {\color{black}
\begin{alignat*}{2}\label{ex10}
	&g(x)= \exp\Bigl(\sum_{i=1}^{d}(-1)^{i+1}x_i \Bigr),
	&&\qquad \widehat{g}(x)= \cos \Bigl(2\pi+2\sum_{i=1}^{d} x_i \Bigr), \\
	&  \widetilde{g}(x)= \prod_{i=0}^{d} \frac{1}{0.9^2+(x_i-0.6)^2}.
\end{alignat*}

Table \ref{tab:13}, \ref{tab:14} and \ref{tab:19} present respectively the computed results of {\bf Test 9} with $d=5,10$ and $N=11,21,41,81,161,321$. } 
It should be noted that the quality of the approximation
also depends on the behavior of the integrand. For very oscillatory and fast growth functions, more 
integration points must be used to achieve good accuracy. In the next section, we shall examine using the regression technique the relationship between the CPU time and the parameter $N$ and the dimension $d$. 
\begin{table}[H]
	\centering
	\begin{tabular}{cccccc}
		\cline{1-6} \noalign{\smallskip}
		\multicolumn{2}{c}{}&\multicolumn{2}{c}{$d=5$ } &\multicolumn{2}{c}{$d=10$ } \\
		\cline{3-6} \noalign{\smallskip}
		\makecell[c]{Mesh size\\ ($h$)}&\makecell[c]{ Points\\ $N$}&\makecell[c]{Relative\\ error} &\makecell[c]{CPU\\time(s)} &\makecell[c]{Relative\\ error} &\makecell[c]{CPU\\time(s)} \\
		\noalign{\smallskip}\hline\noalign{\smallskip}
		0.1  & $11$ & $2.7744\times 10^{-6}$   & 0.1756857   & $5.5489\times 10^{-6}$   & 0.4397261\\
		0.05  & $21$ & $1.7355\times 10^{-7}$   &  0.3406698   & $3.4711\times 10^{-7}$   & 1.0226371\\
		0.025   &$41$  & $1.0849\times 10^{-8}$ &  0.7929361   & $2.1699\times 10^{-8}$ &  4.0287093  \\
		0.0125   & $81$& $6.7815\times 10^{-10}$&  2.1783682& $1.3563\times 10^{-9}$& 31.1292788  \\	
		0.00625  &$161$ &  $4.2384\times 10^{-11}$ & 31.7971237  &  $8.4768\times 10^{-11}$ & 141.9877281\\
		0.003125  &$321$ &  $2.6479\times 10^{-12}$ & 136.4977085   &  $5.2993\times 10^{-12}$ & 550.772326  \\		
		\noalign{\smallskip}\hline
	\end{tabular}
	\caption{Performance comparison of the MDI algorithm with $d=5,10$ and $N=11,21,41,81,161,321$ for approximating $I(g,\Omega)$.}   \label{tab:13}  
\end{table}

\begin{table}[H]
	\centering
	\begin{tabular}{cccccc}
		\cline{1-6} \noalign{\smallskip}
		\multicolumn{2}{c}{}&\multicolumn{2}{c}{$d=5$ } &\multicolumn{2}{c}{$d=10$ } \\
		\cline{3-6} \noalign{\smallskip}
		\makecell[c]{Mesh size\\ ($h$)}&\makecell[c]{ Points\\ $N$}&\makecell[c]{Relative\\ error} &\makecell[c]{CPU\\time(s)} &\makecell[c]{Relative\\ error} &\makecell[c]{CPU\\time(s)} \\
		\noalign{\smallskip}\hline\noalign{\smallskip}
		0.1  & $11$ &$4.4657\times 10^{-5}$ &  0.1832662   &$8.9317\times 10^{-5}$ & 0.4046452\\
		0.05  & $21$ &$2.7810\times 10^{-6}$ & 0.3478922   &$5.5621\times 10^{-6}$ &1.0245456 \\
		0.025   &$41$   &$1.7366\times 10^{-7}$   & 0.7666344  &$3.4732\times 10^{-7}$   &  4.5341739  \\
		0.0125   & $81$&$1.0851\times 10^{-8}$    & 2.8629733 &$2.1703\times 10^{-8}$    & 35.0469255  \\	
		0.00625  &$161$  & $6.7818\times 10^{-10}$ &26.3679753  & $ 1.3563\times 10^{-9}$ &150.875945 \\
		0.003125  &$321$  & $4.2383\times 10^{-11}$ & 146.441442   & $8.4768\times 10^{-11}$ &568.693914 \\		
		\noalign{\smallskip}\hline
	\end{tabular}
\caption{Performance comparison of the MDI algorithm with $d=5,10$ and $N=11,21,41,81,161,321$ for approximating $I(\widehat{g},\Omega)$.}\label{tab:14}    
\end{table}

\begin{table}[H]
	\centering
	\begin{tabular}{cccccc}
		\cline{1-6} \noalign{\smallskip}
		\multicolumn{2}{c}{}&\multicolumn{2}{c}{$d=5$ } &\multicolumn{2}{c}{$d=10$ } \\
		\cline{3-6} \noalign{\smallskip}
		\makecell[c]{Mesh size\\ ($h$)}&\makecell[c]{ Points\\ $N$}&\makecell[c]{Relative\\ error} &\makecell[c]{CPU\\time(s)} &\makecell[c]{Relative\\ error} &\makecell[c]{CPU\\time(s)} \\
		\noalign{\smallskip}\hline\noalign{\smallskip}
		0.1  & $11$ & $2.6251\times 10^{-5}$   &  0.1669653   &$5.2504\times 10^{-5}$ & 0.273292\\
		0.05  & $21$ & $1.6339\times 10^{-6}$   & 0.2529685   &$3.2679\times 10^{-6}$ &0.4857094 \\
		0.025   &$41$  & $1.0200\times 10^{-7}$ &  0.4683767   &$2.0401\times 10^{-7}$   & 0.8593056  \\
		0.0125   & $81$& $6.3736\times 10^{-9}$& 0.8346902&$1.2747\times 10^{-8}$    & 1.7529975  \\	
		0.00625  &$161$ &  $3.9832\times 10^{-10}$ & 1.6505776  & $ 7.9664\times 10^{-10}$ &3.6539963 \\
		0.003125  &$321$ &  $ 2.4894\times 10^{-11}$ & 3.6465664   & $ 4.9788\times 10^{-11}$ &9.5368530 \\		
		\noalign{\smallskip}\hline
	\end{tabular}
\caption{Performance comparison of the MDI algorithm with $d=5,10$ and $N=11,21,41,81,161,321$ for approximating $I(\widetilde{g},\Omega)$.}\label{tab:19} 
\end{table}

\section{Computational complexity}\label{sec-5}
{\color{black} 
\subsection{The relationship between the CPU time and $N$} \label{sec-5.1}
In this subsection, we examine  the relationship between the CPU time and parameter $N$ using the regression technique
based on the test data.
\begin{figure}[H]
	\centerline{
		\includegraphics[width=1.8in,height=1.7in]{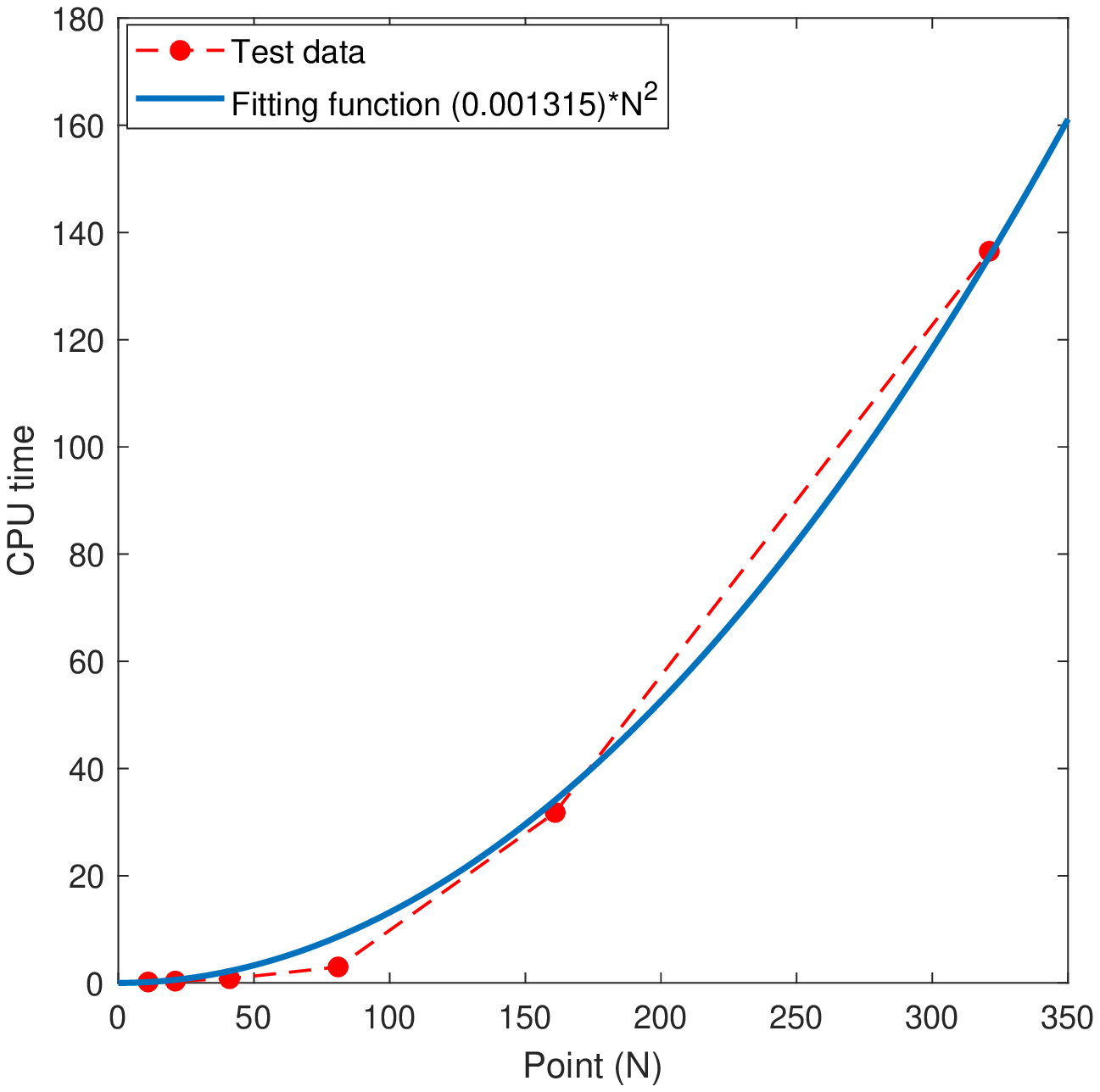}
		\includegraphics[width=1.8in,height=1.7in]{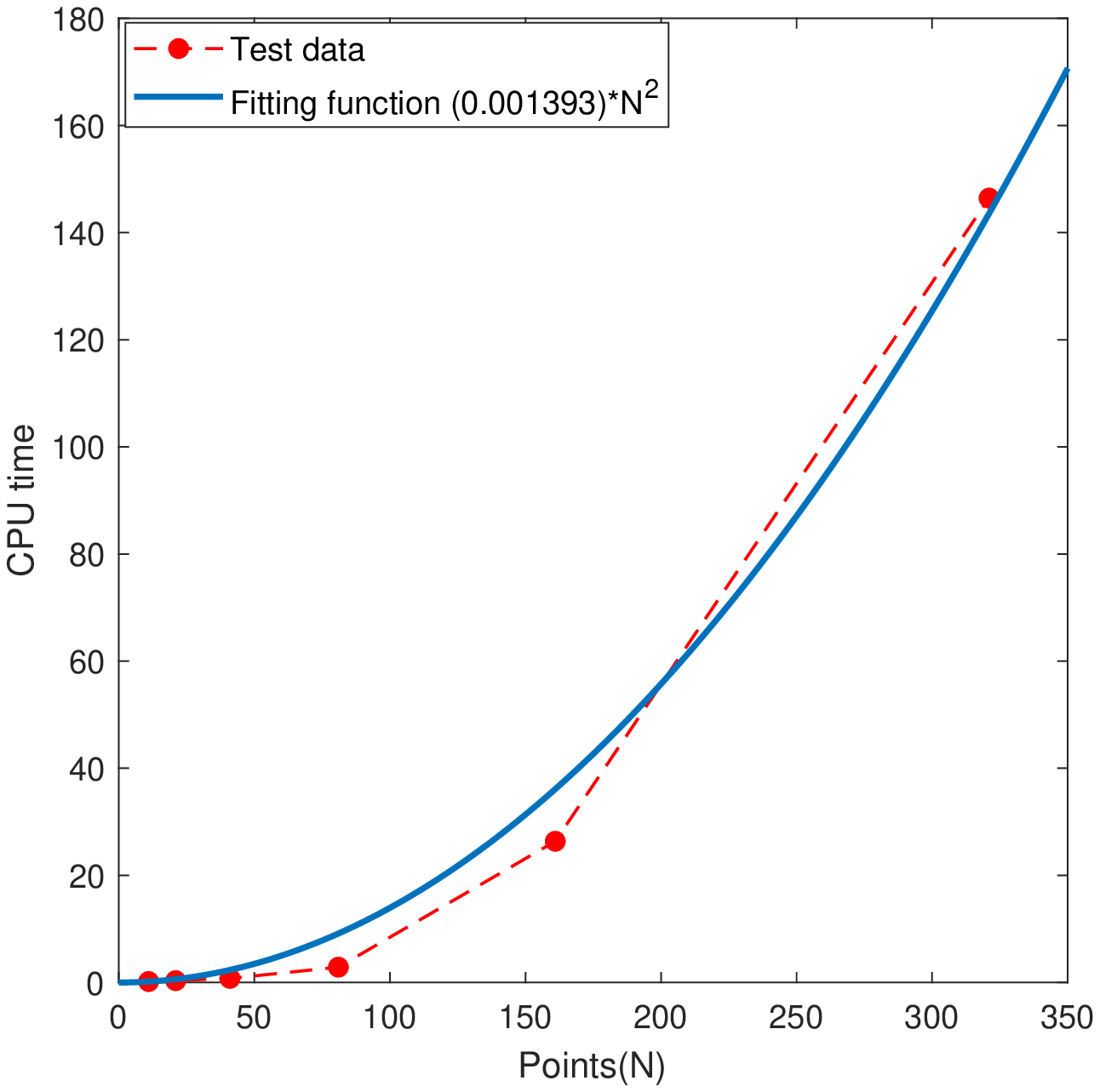}
		\includegraphics[width=1.8in,height=1.7in]{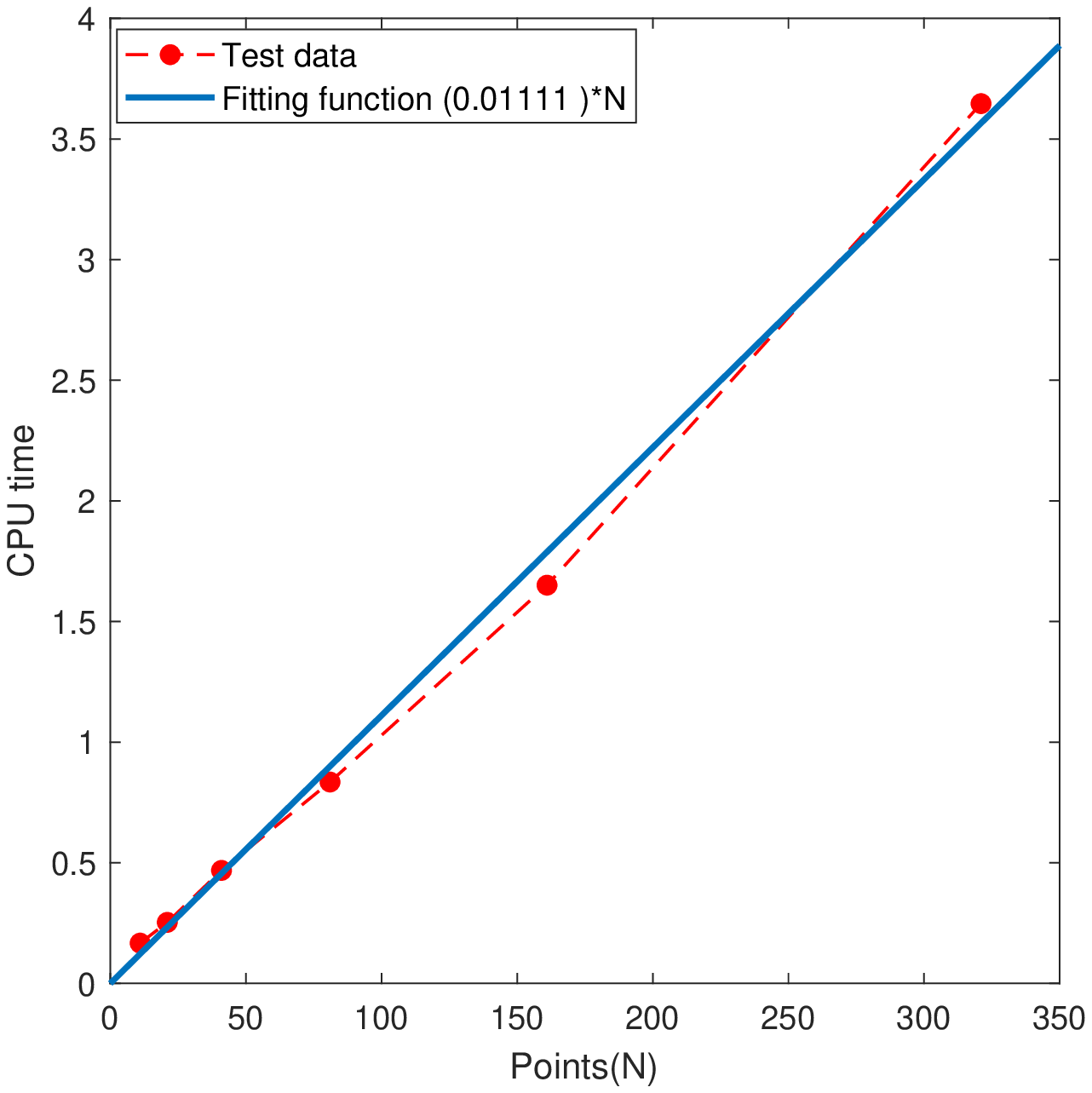}
	}
	\caption{The relationship between the CPU time and parameter $N$ when $d=5$: $I(g,\Omega)$ (left), $I(\widehat{g},\Omega)$ (middle), $I(\widetilde{g},\Omega)$ (right).}
	\label{fig.2}       
\end{figure}
\begin{figure}[H]
	\centerline{
		\includegraphics[width=1.8in,height=1.7in]{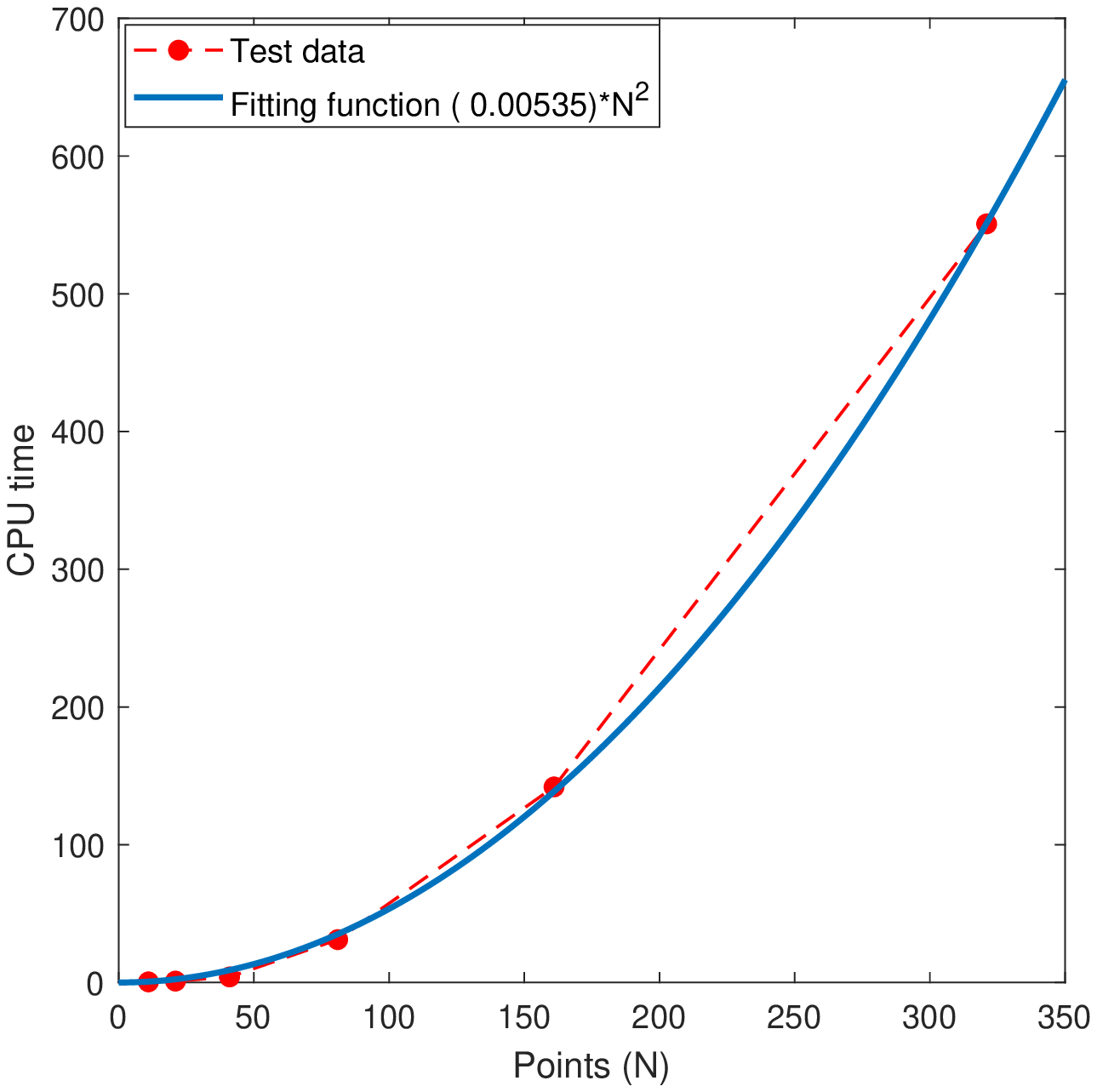}
		\includegraphics[width=1.8in,height=1.7in]{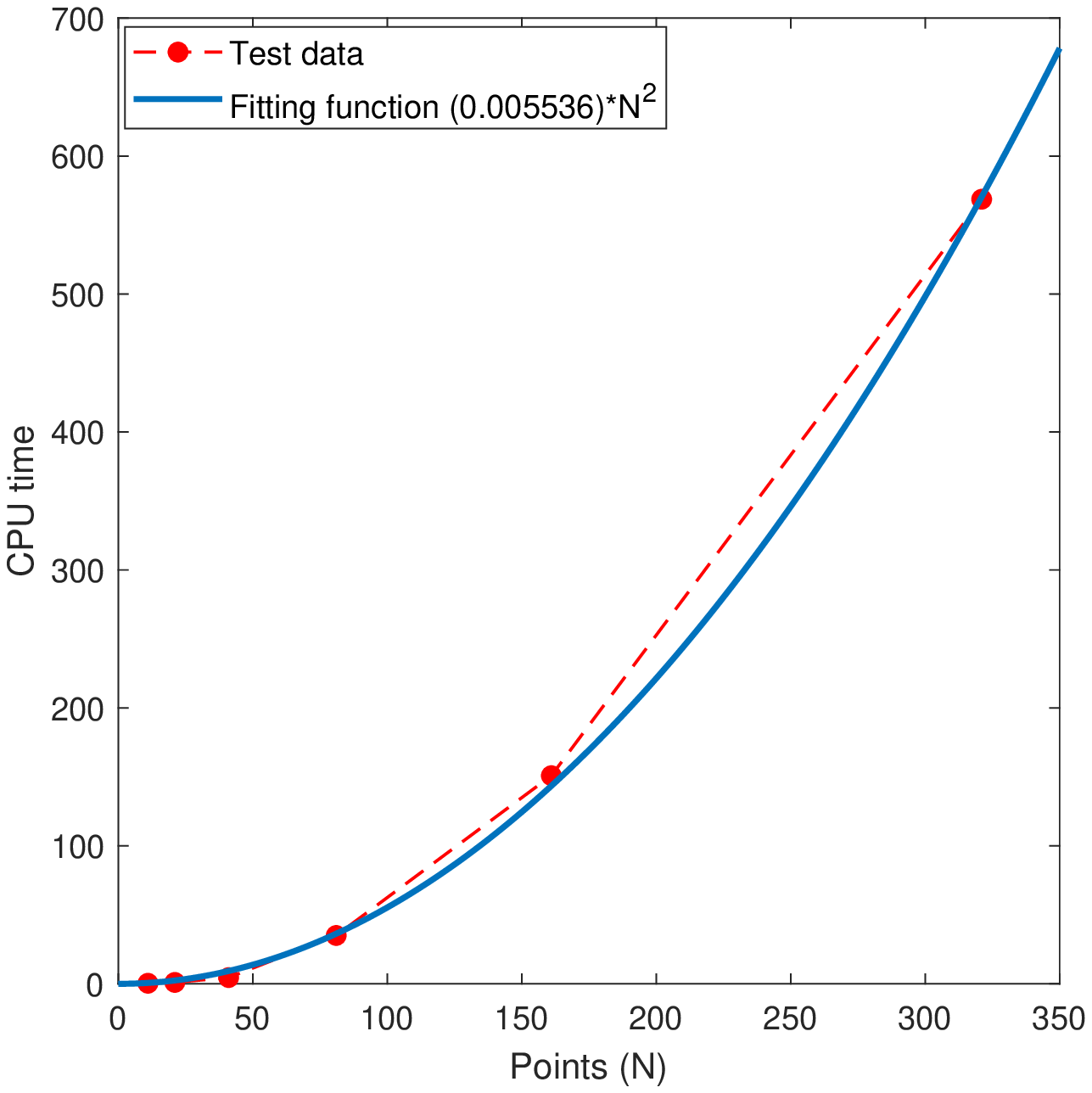}
		\includegraphics[width=1.8in,height=1.7in]{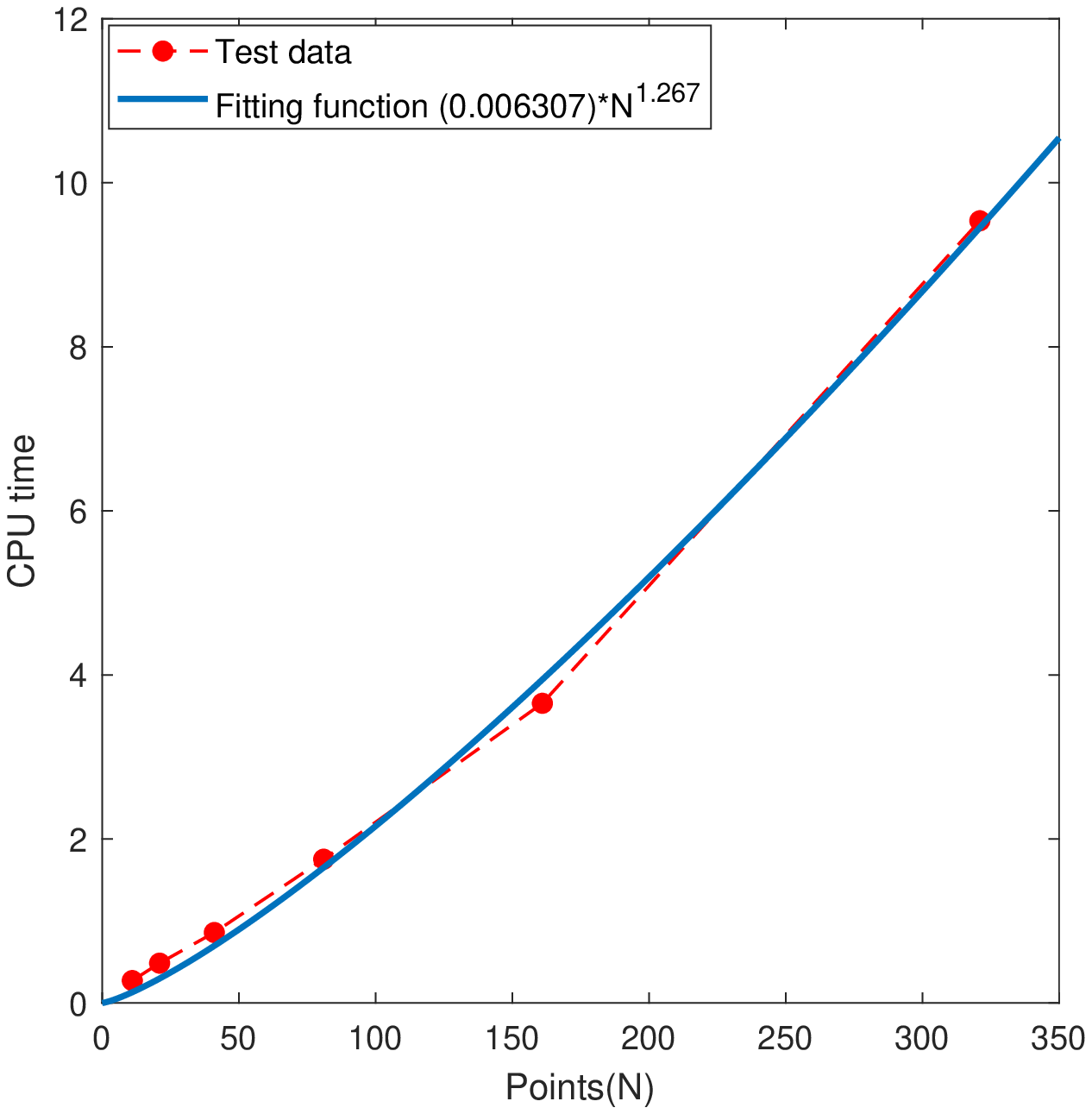}
	}
	\caption{The relationship between the CPU time and parameter $N$ when $d=10$: $I(g,\Omega)$ (left), $I(\widehat{g},\Omega)$ (middle), $I(\widetilde{g},\Omega)$ (right).}
	\label{fig.3}       
\end{figure}
\begin{table}[H]
	\centering
	\begin{tabular}{llllll}
		\toprule
		\makecell[c]{Integrand\\}&\makecell[c]{$r$} &\makecell[c]{$m$} &\makecell[c]{$d$} &\makecell[c]{Fitting function}  &\makecell[c]{R-square\\}\\
		
		\midrule
		$g(x)$ &2&  1& 5 & $h_1(N) =(0.001315)*N^2$ &0.9973   \\
		$\widehat{g}(x)$&2&  1& 5  & $h_2(N)=(0.001393)*N^2$ &0.9914    \\
		$\widetilde{g}(x)$&2&  1& 5  & $h_3(N)=(0.01111 )*N$ &0.9963     \\ 
		$g(x)$&2&  1& 10  &$h_4(N) =( 0.00535)*N^2$  &0.9998    \\
		$\widehat{g}(x)$&2&  1& 10  &$h_5(N) =(0.005536)*N^2$  &0.9997    \\
		$\widetilde{g}(x)$&2&  1& 10  &$h_6(N)=(0.006307)*N^{1.267}$  &0.9971    \\
		\bottomrule
	\end{tabular}
	\caption{The relationship between the CPU time and parameter $N$.} \label{tab:20}
\end{table}

Fig. \ref{fig.2} and \ref{fig.3} show the CPU time as a function of $N$ obtained by the least square method and the fitting functions are given in Table \ref{tab:20}. All the results indicate that the CPU time grows at most quadratically in $N$.

}

\subsection{ the relationship between the CPU time and the dimension $d$ }
\label{sec-5.2} 
 
Recall that the computational complexity of tensor product methods is of the exponential order $O(N^d)$. 
The numerical tests presented above overwhelmingly and consistently show that the MDI algorithm 
has hidden capability to overcome the curse of dimensionality faced by tensor product methods. 
The goal of the next test is to find out the computational complexity (in terms of CPU time as a function
of $d$) using the least square method based on numerical test data.

\medskip
{\bf Test 10.} Let $\Omega=[0,1]^d$, we consider the following five integrands: 
\begin{alignat*}{2}\label{ex11}
&g_1(x)= \exp\Bigl(\sum_{i=1}^{d}(-1)^{i+1}x_i \Bigr),
&&\qquad g_2(x)= \prod_{i=0}^{d} \frac{1}{0.9^2+(x_i-0.6)^2},\\
&g_3(x)= \frac{1}{\sqrt{2\pi}}\exp\Bigl(-\frac{1}{2}|x|^2\Bigr), 
&&\qquad g_4(x)= \cos \Bigl(2\pi+2\sum_{i=1}^{d} x_i \Bigr), \\
&g_5(x)= \exp\Bigl(5\sum_{i=1}^{d}(-1)^{i+1}x_i^2 \Bigr).
\end{alignat*}

 	Fig. \ref{fig.11} displays the the CPU time as functions of $d$ obtained by the least square method whose analytical expressions are given in Table \ref{tab:15}. We note that the parameters of the MDI algorithm only affect the coefficients of the fitting functions, but not the order.

\begin{table}[H]
	\centering
	\begin{tabular}{lccclc}
	    \cline{1-6} \noalign{\smallskip}
		\makecell[c]{Integrand\\}&\makecell[c]{$r$} &\makecell[c]{$m$} &\makecell[c]{$N$} &\makecell[c]{Fitting function}  &\makecell[c]{R-square\\}\\
		\noalign{\smallskip}\hline\noalign{\smallskip}
	
		$g_1$
	 	  &1&  1& 11  & $f_1=(1.903e-06)*N^2d^3$ &0.9966   \\
	 	  &2&  1& 7  & $f_2=(1.219e-06)*N^2d^3$ &0.9961   \\
		  &2&  1& 11  & $f_2=(1.761e-06)*N^2d^3$ &0.9964   \\
		   &3&  1& 3  & $f_3=(1.01e-06)*N^2d^3$ &0.9978     \\ 
		  &4&  1& 10  &$f_4=(1.506e-06)*N^2d^3$  &0.9947    \\
		  \hline
	
		$g_2$
		 &1&  1& 11  & $f_5=0.0002787*N^2d^1$ &0.9922     \\
		 &2&  1& 7  & $f_6=(6.154e-05)*N^2d^{1.358}$ &0.9991    \\
		&2&  1& 11  & $f_6=0.0002779*N^2d^1$ &0.9898    \\
		&2&  2& 11  & $f_9=0.001737*N^2d^1$ &0.9937     \\
		&2&  3& 11  & $f_{10}=0.01517*N^2d^1$ &0.9808     \\
		&3&  1& 10  & $f_7=(2.629e-07)*N^2d^3$ &0.9932     \\
		&4&  1& 10  & $f_8=0.000271*N^2d^1$ &0.9952     \\
	    \hline
	    
		$g_3$ 
		&1&  1& 11  &  $f_{11}=(2.629e-05)*N^2d^3$ &0.9977    \\
		 &2&  1& 11  & $f_{12}=(2.682e-05)*N^2d^3$ &0.9995     \\
		&2&  2& 11  & $f_{13}=0.0001357*N^2d^3$ &0.9929     \\
		&2&  3& 11  & $f_{14}=0.0008642*N^2d^3$ &0.9974     \\
		&4&  1& 10  & $f_{15}=(1.123e-05)*N^2d^3$ &0.9903     \\
		\hline
		
		$g_4$  
		
		&2&  1& 11  & $f_{16}=(1.633e-06)*N^2d^3$ &0.9990      \\
		&2&  1& 21  & $f_{17}=(1.74e-06)*N^2d^3$ &0.9966      \\
		&3&  1& 3  & $f_{18}=(7.904e-07 )*N^2d^3$ &1.0000      \\
		\hline
		$g_5$ 
		 &2&  1& 11  &$f_{19}=(2.136e-05)*N^2d^3$ &0.9968       \\
		  &2&  1& 21  & $f_{20}=(4.791e-05)*N^2d^3$ &0.9994      \\
		\noalign{\smallskip}\hline
	\end{tabular}
\caption{The relationship between CPU time and the integral dimension $d$.}\label{tab:15}     
\end{table}
	{\color{black}we quantitatively characterize the performance of the fitted curve by the $R$-square in Matlab, which is defined as $R$-square$=1-\frac{\sum_{i}^n(y_i-\widehat{y}_i)^2}{\sum_{i}^n(y_i-\overline{y})^2}$. Where $y_i$ represents the test data, $\widehat{y}_i$ refers to the predicted value, and $\overline{y}$ indicates the mean value of $y_i$.
	Table \ref{tab:15} also shows that the R-square of all fitting functions is very close to $1$, which indicates the 
	fitting function is a quite accurate.} These results indicate that the CPU time grows at most cubically in $d$. Combining the results of {\bf Test 9} in Section \ref{sec-4.3} we conclude that the CPU time required by the proposed MDI algorithm grows at most in the polynomial order $O(d^3N^2)$. 

\begin{figure}[H]
	\centerline{
	 \includegraphics[width=1.75in,height=1.65in]{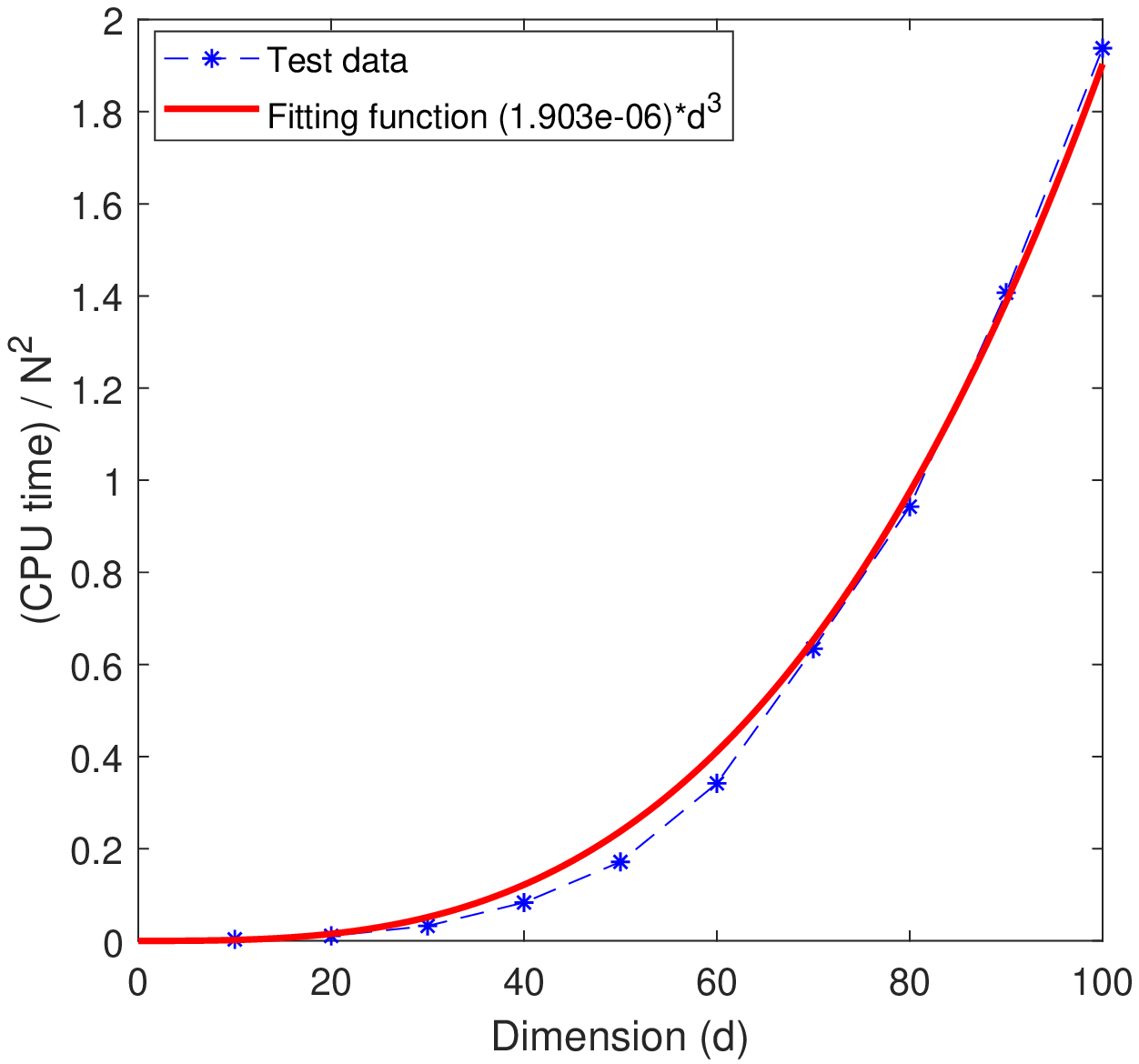}
	 \includegraphics[width=1.75in,height=1.65in]{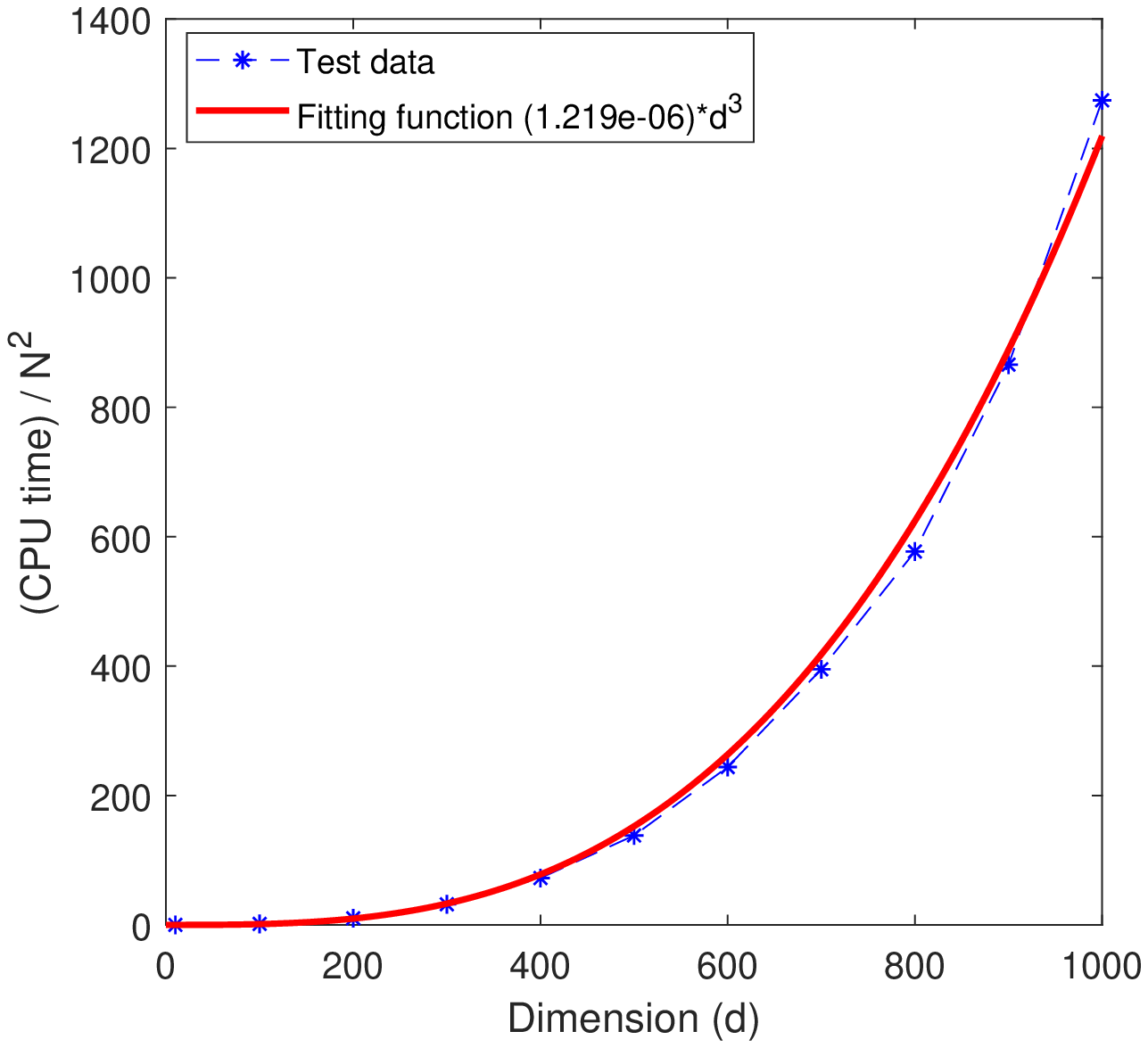}
	 \includegraphics[width=1.75in,height=1.65in]{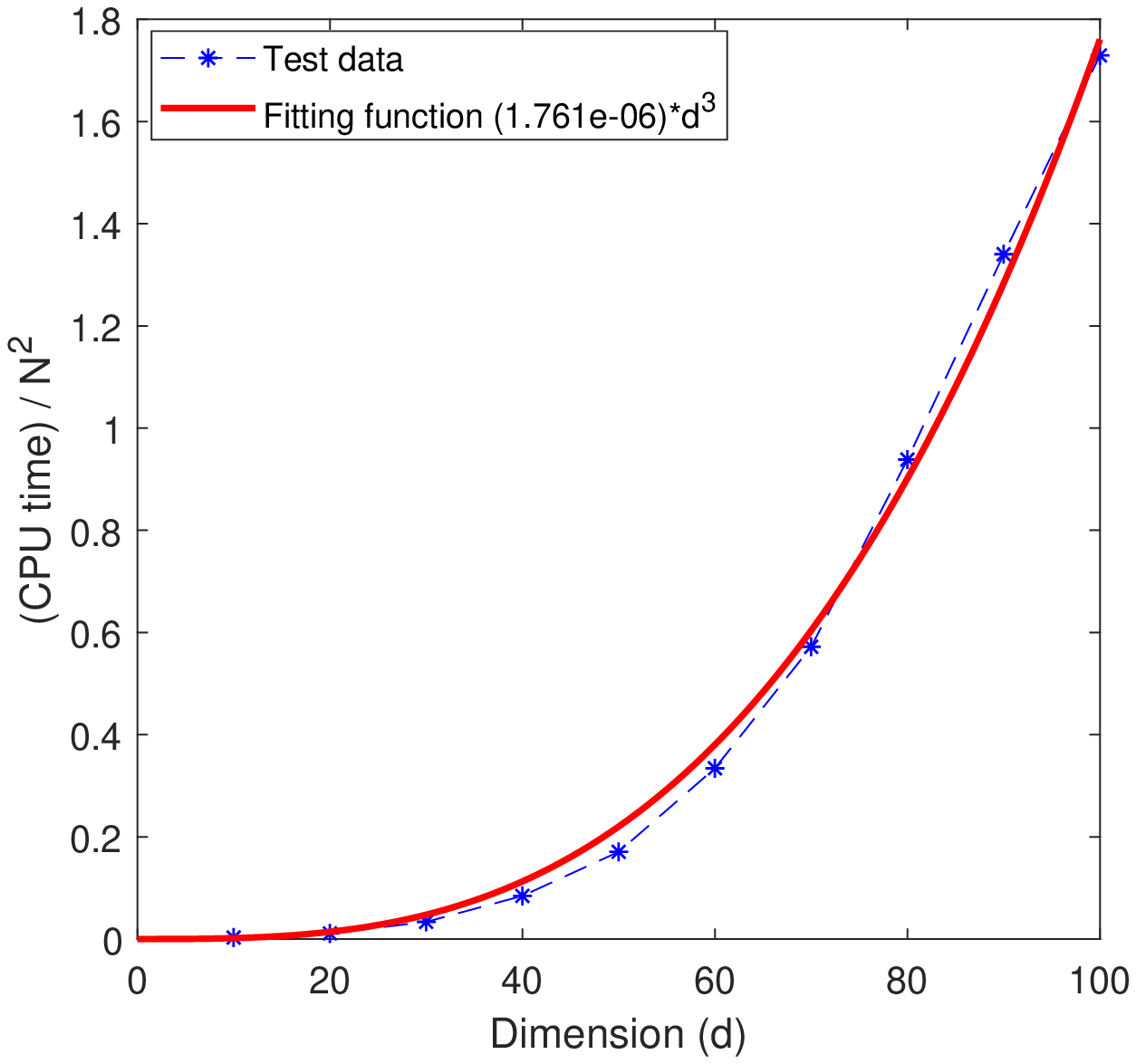}
    }
      \centerline{
      	\includegraphics[width=1.75in,height=1.65in]{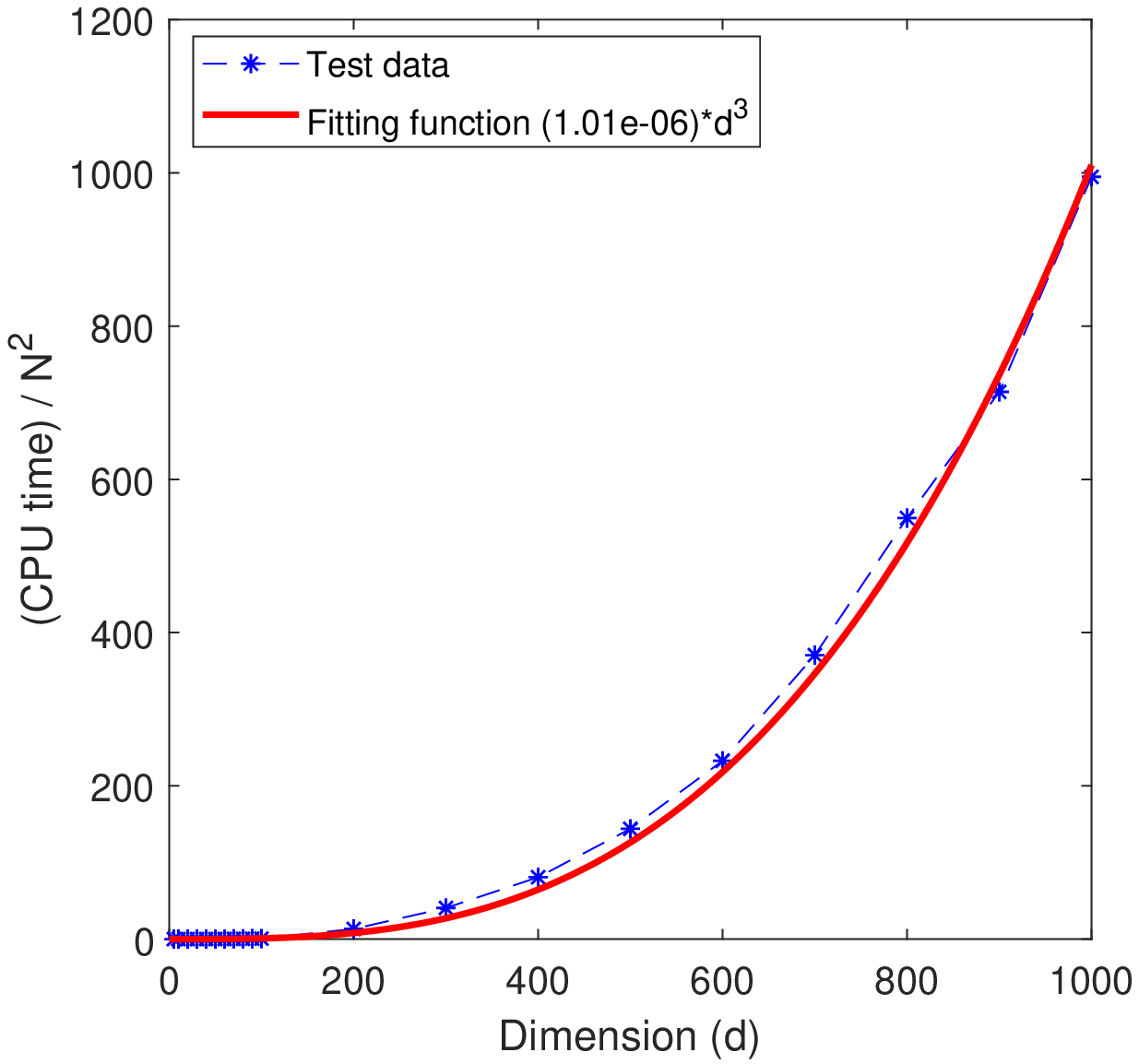}
    	\includegraphics[width=1.75in,height=1.65in]{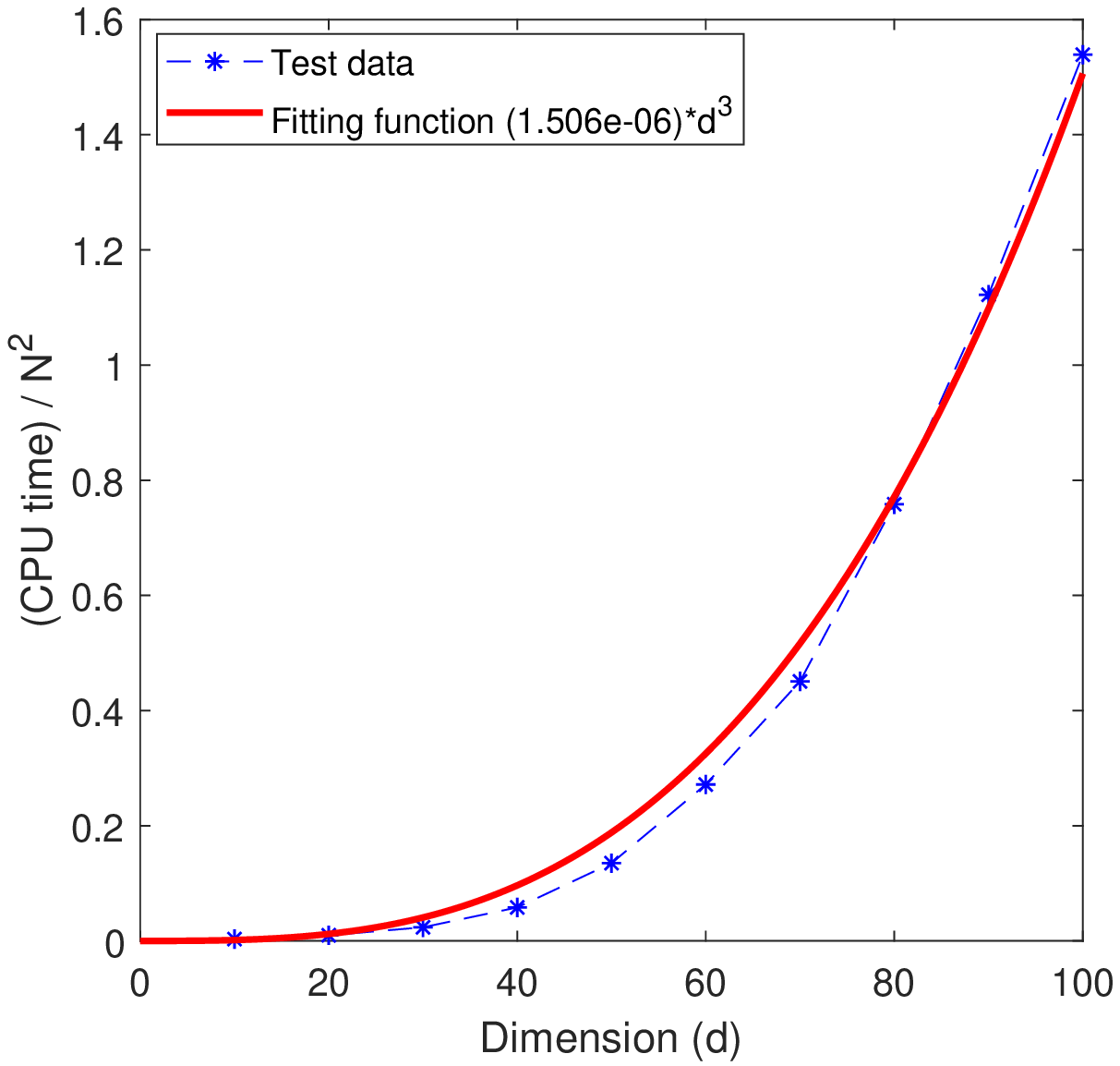}
    	\includegraphics[width=1.75in,height=1.65in]{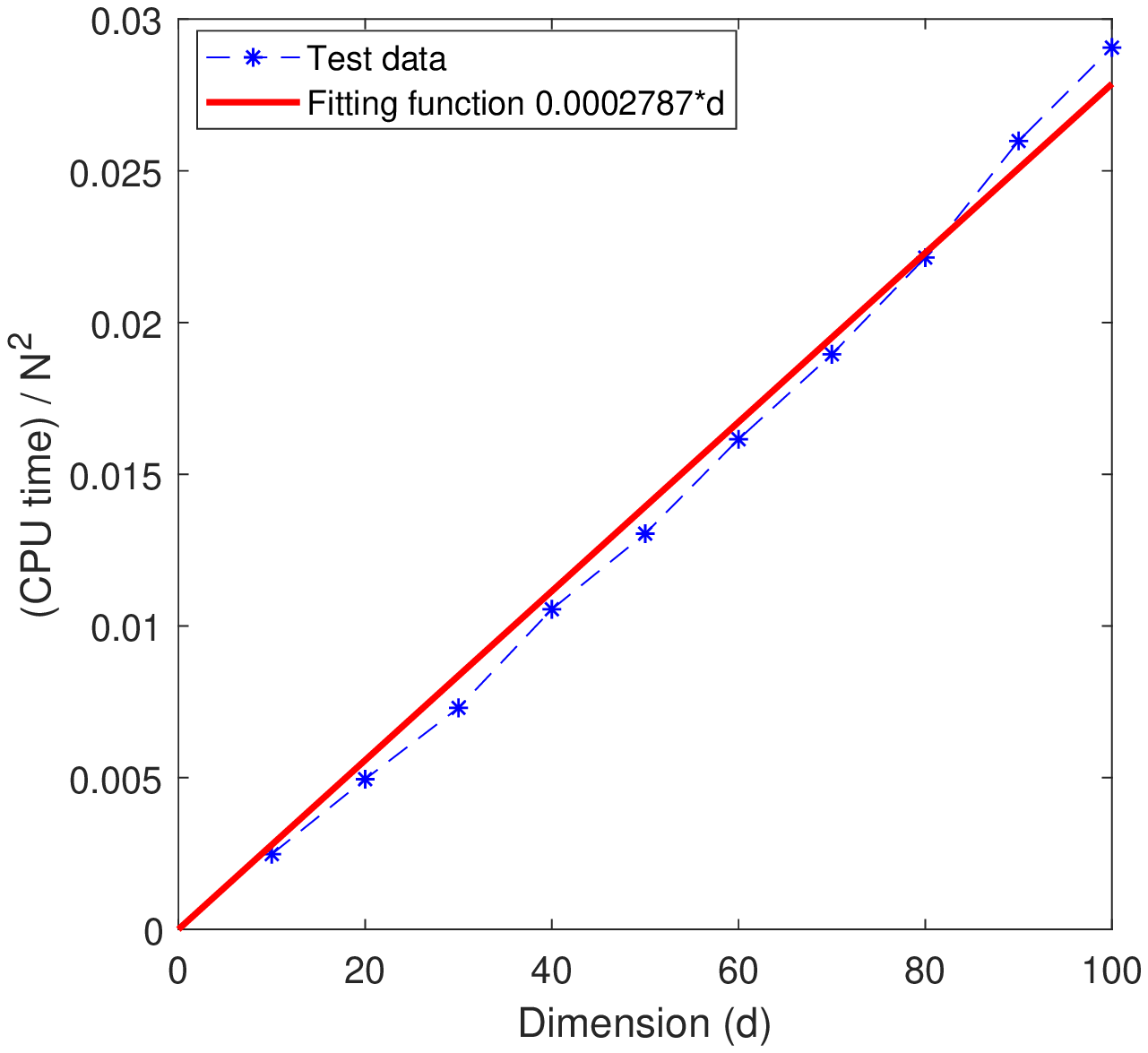}
    }

	\centerline{
	\includegraphics[width=1.75in,height=1.65in]{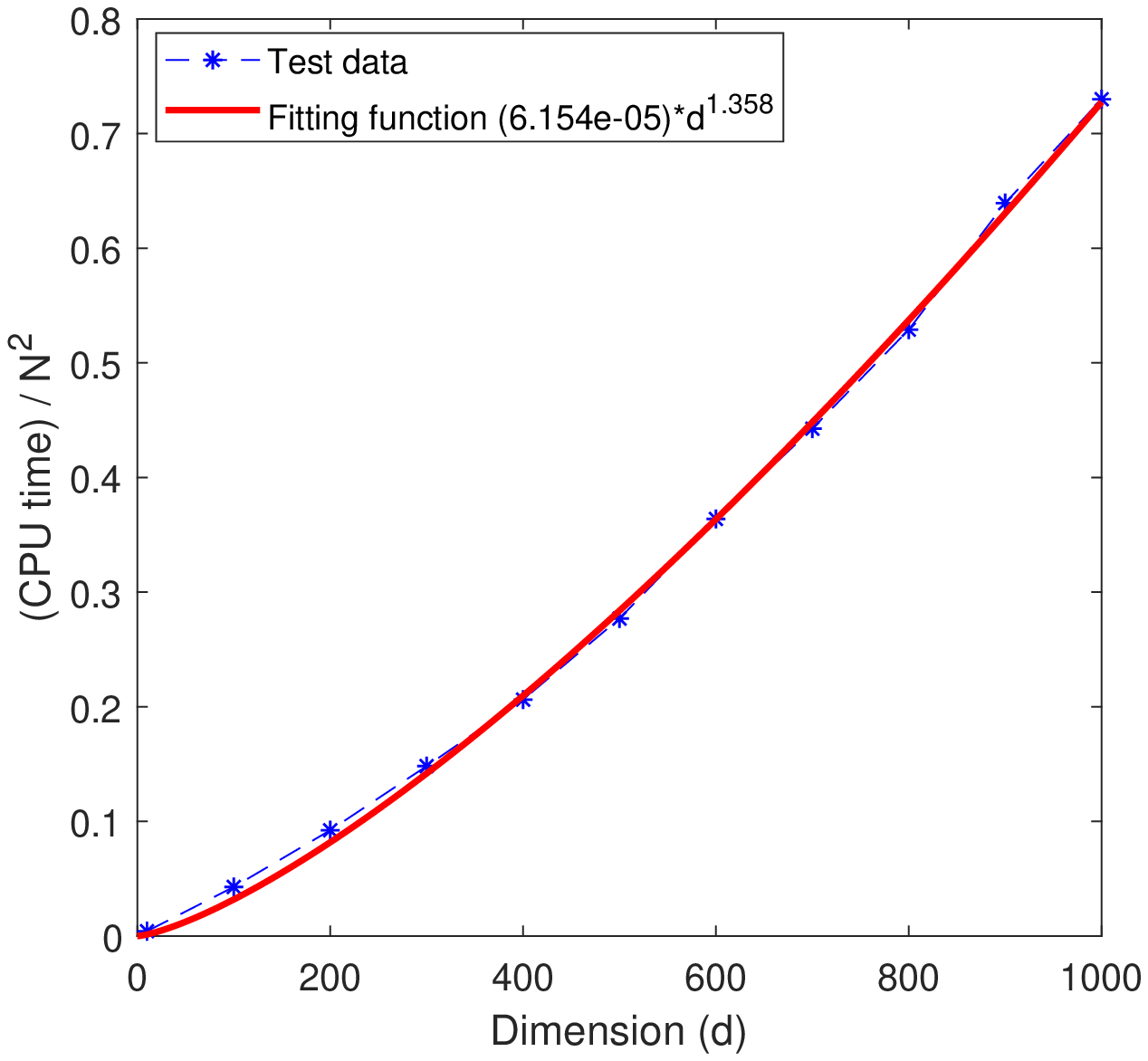}
	\includegraphics[width=1.75in,height=1.65in]{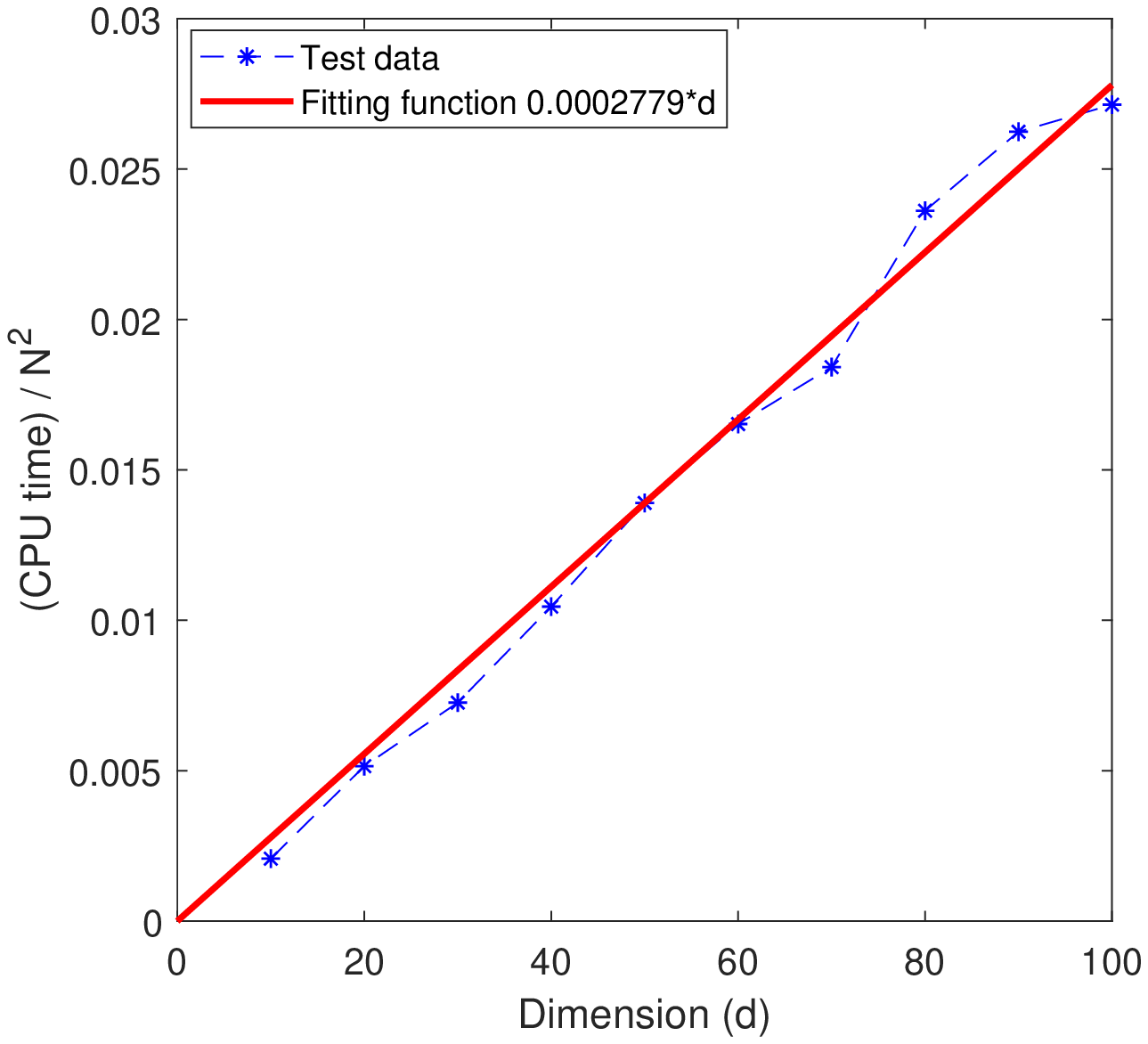}
	\includegraphics[width=1.75in,height=1.65in]{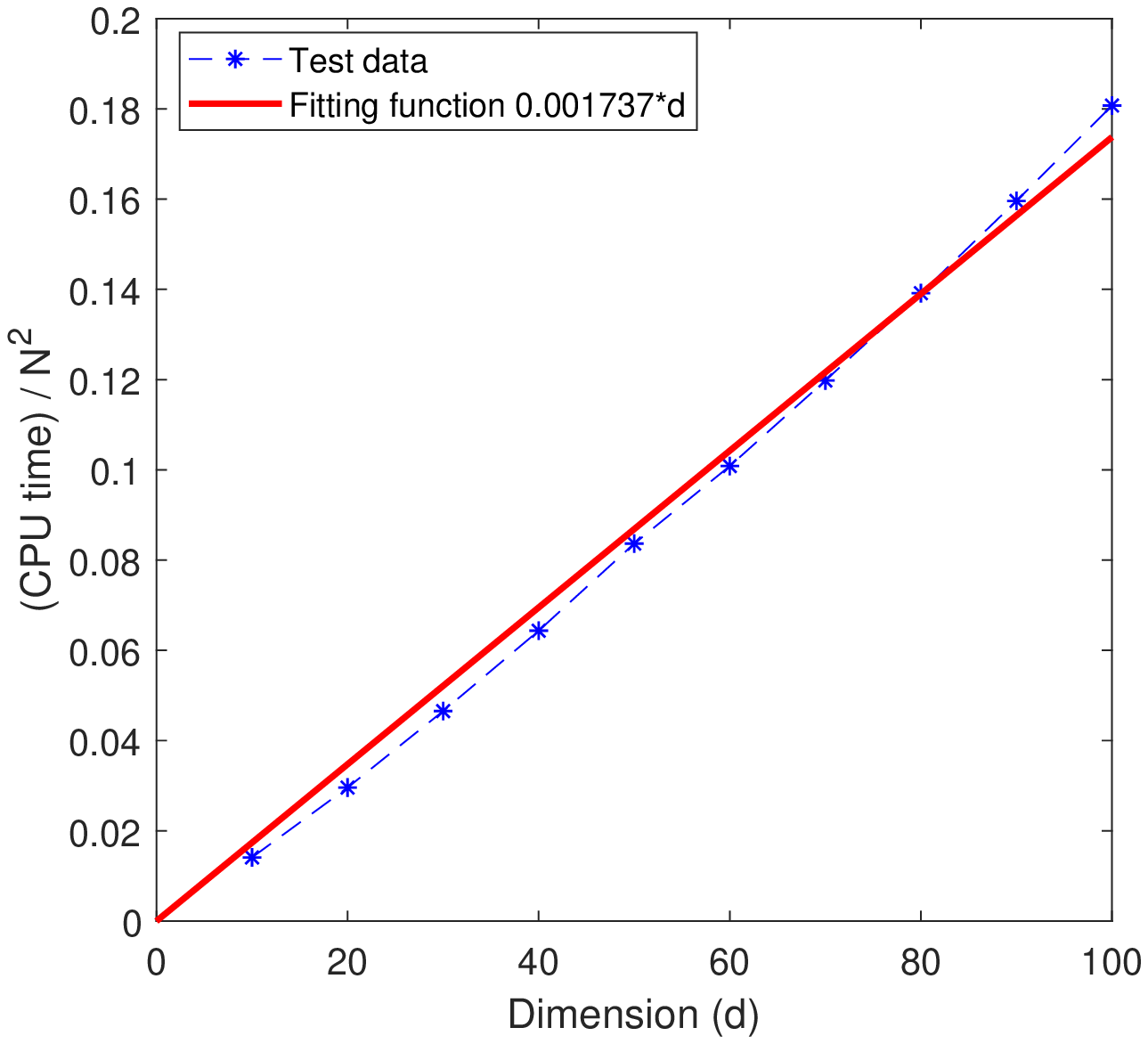}
}
\centerline{
	\includegraphics[width=1.75in,height=1.65in]{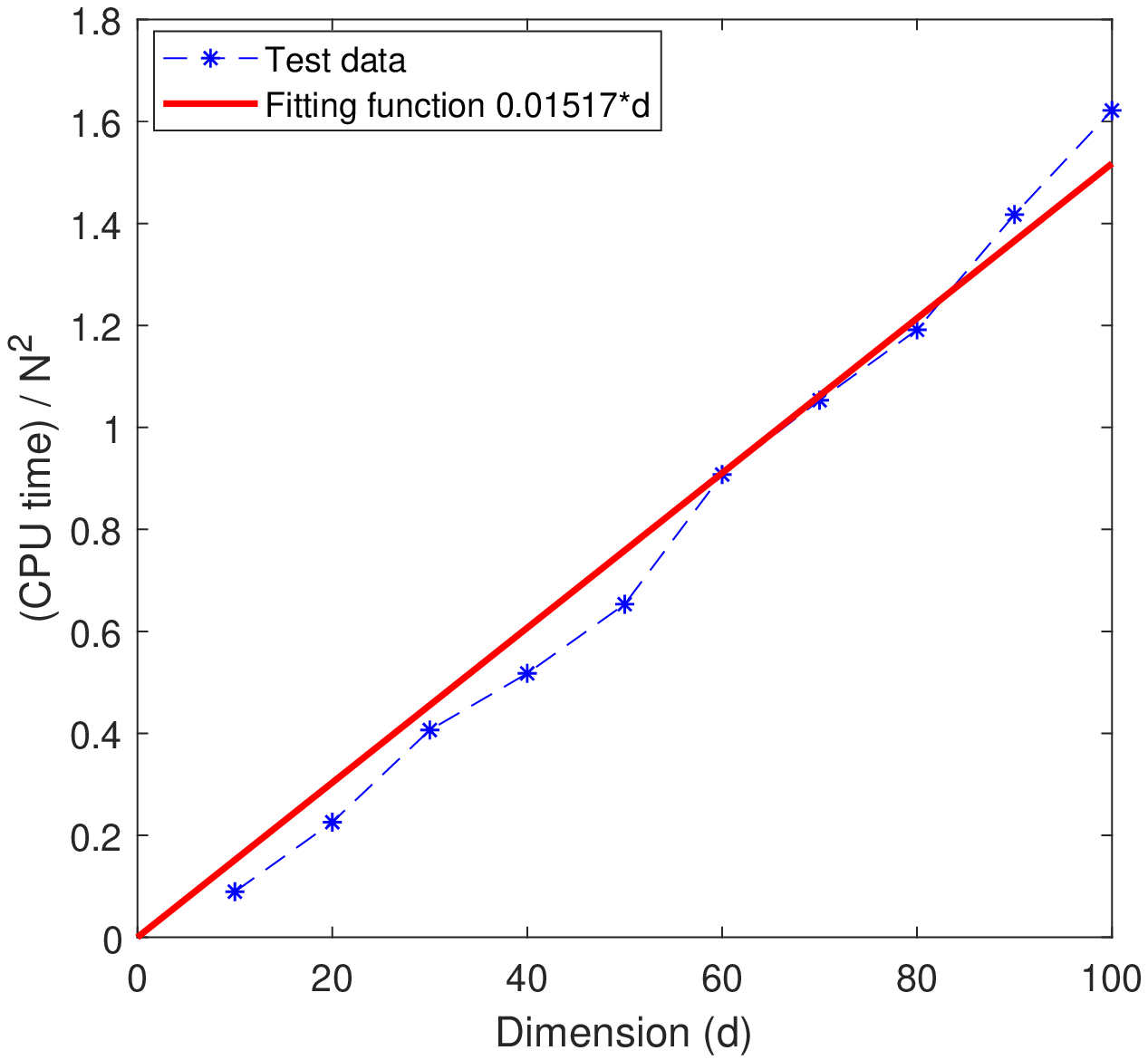}
    \includegraphics[width=1.75in,height=1.65in]{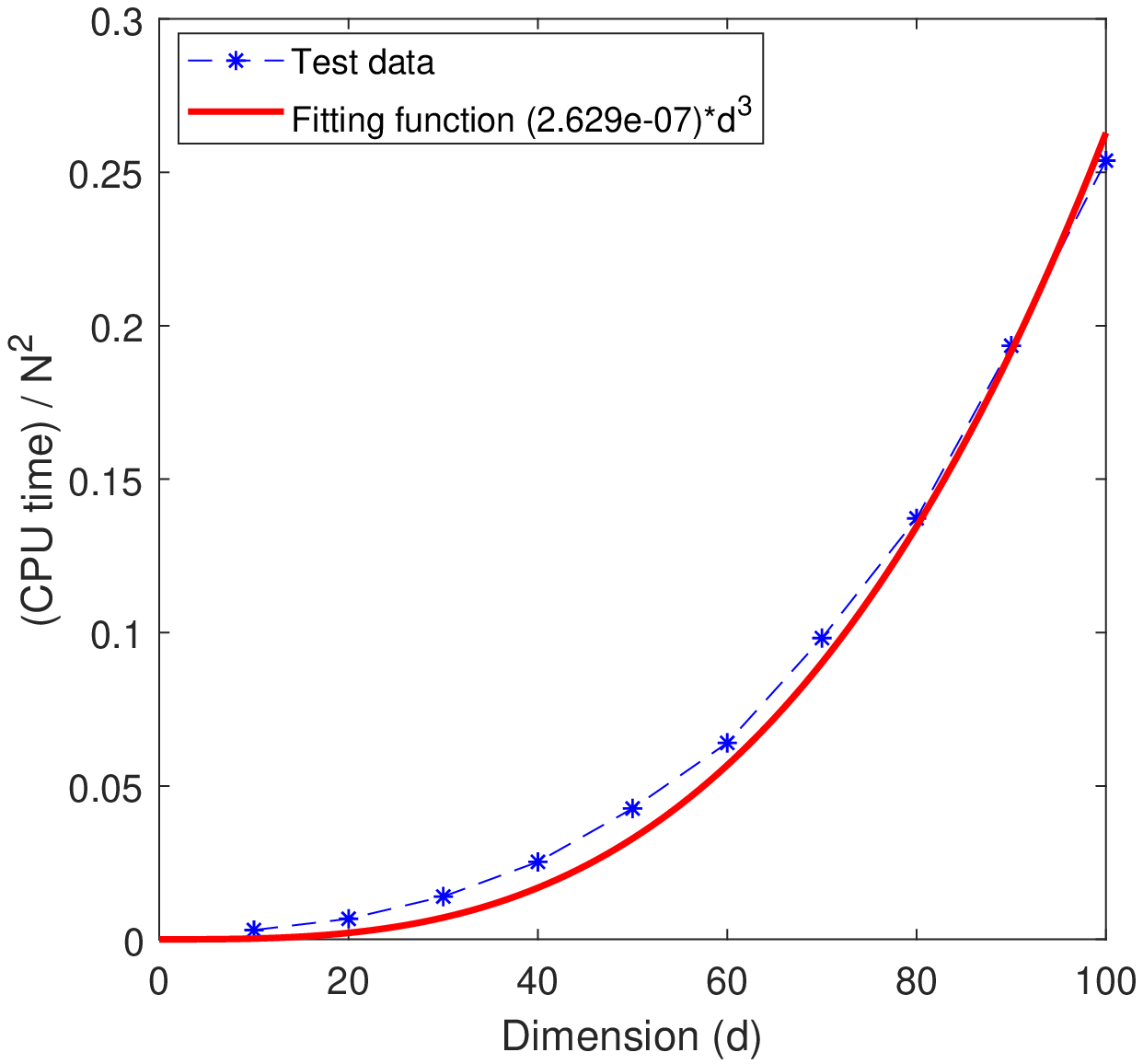}
	\includegraphics[width=1.75in,height=1.65in]{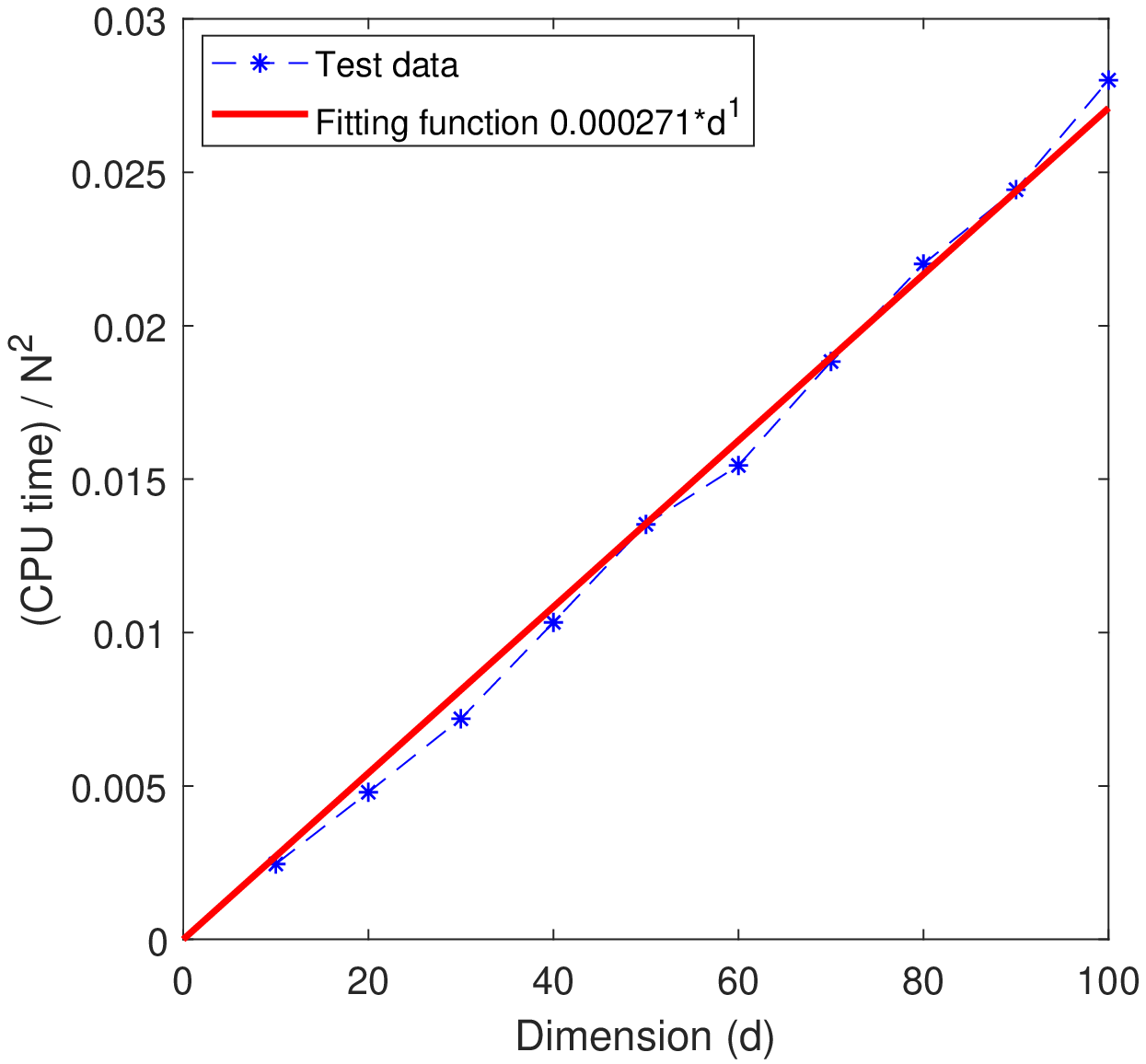}
} 
\end{figure}

\begin{figure}[H]      
		\centerline{
		\includegraphics[width=1.8in,height=1.7in]{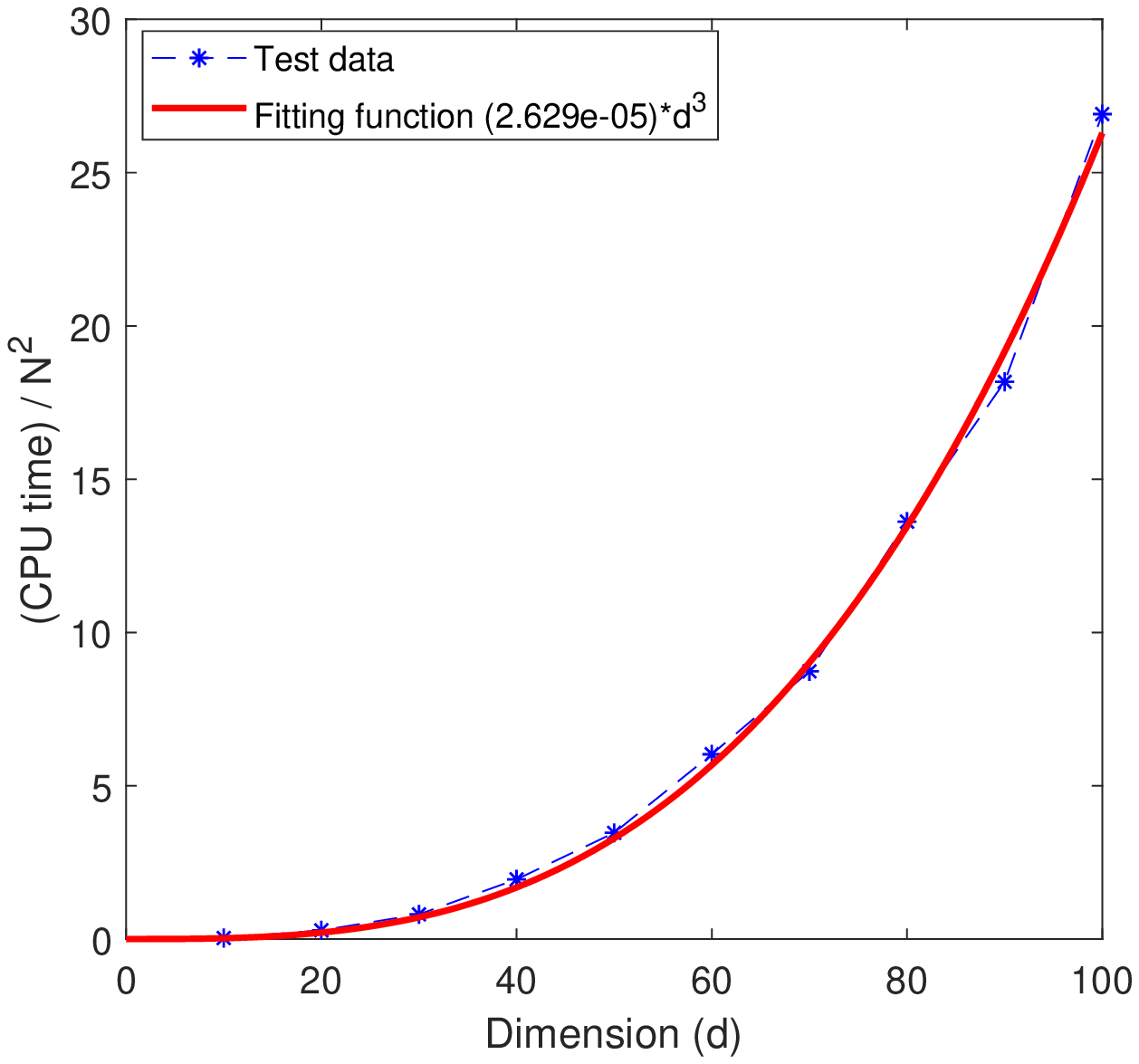}
		\includegraphics[width=1.8in,height=1.7in]{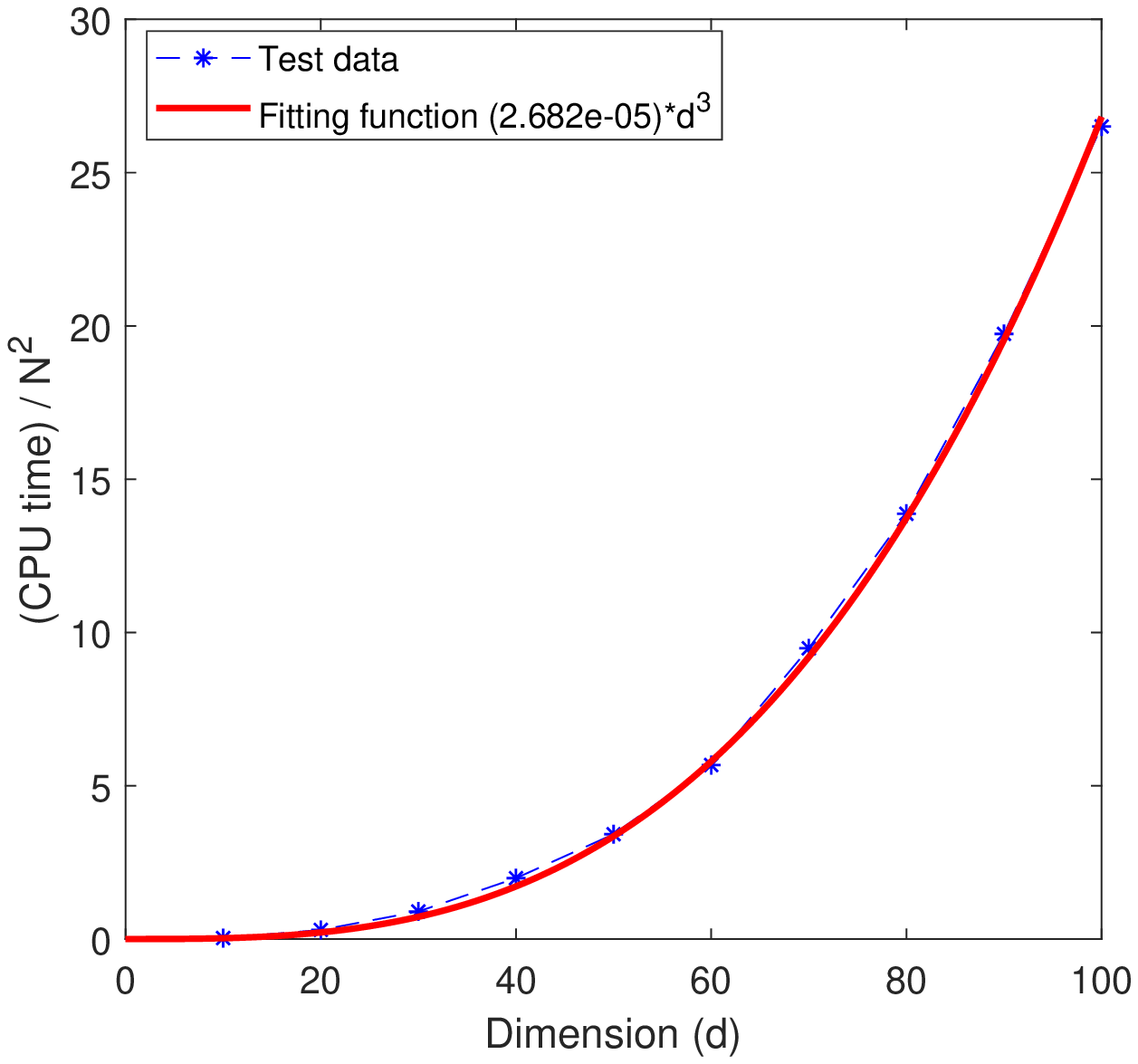}
		\includegraphics[width=1.8in,height=1.7in]{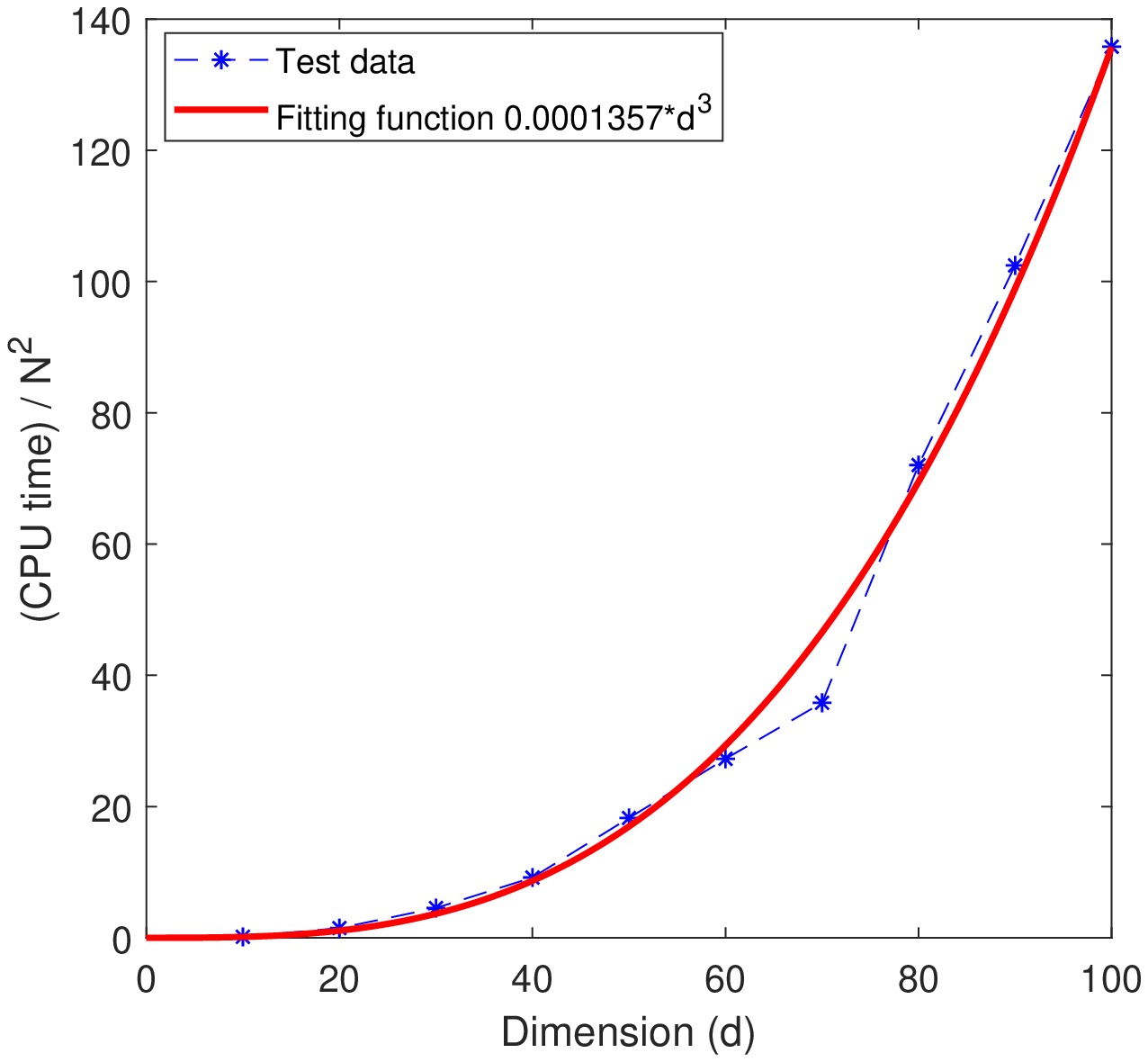}
	}
	\centerline{
    	\includegraphics[width=1.8in,height=1.7in]{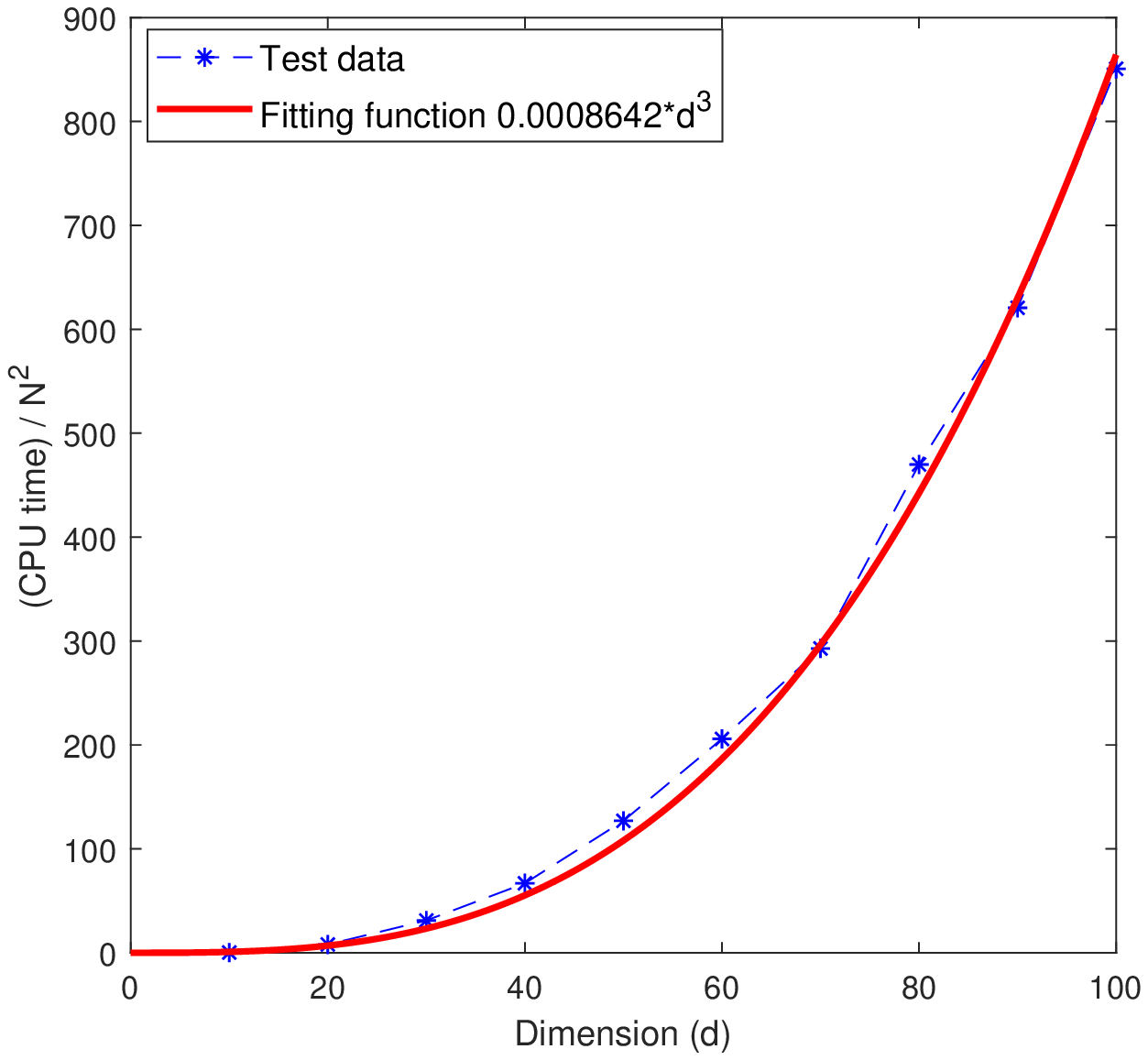}
        \includegraphics[width=1.8in,height=1.7in]{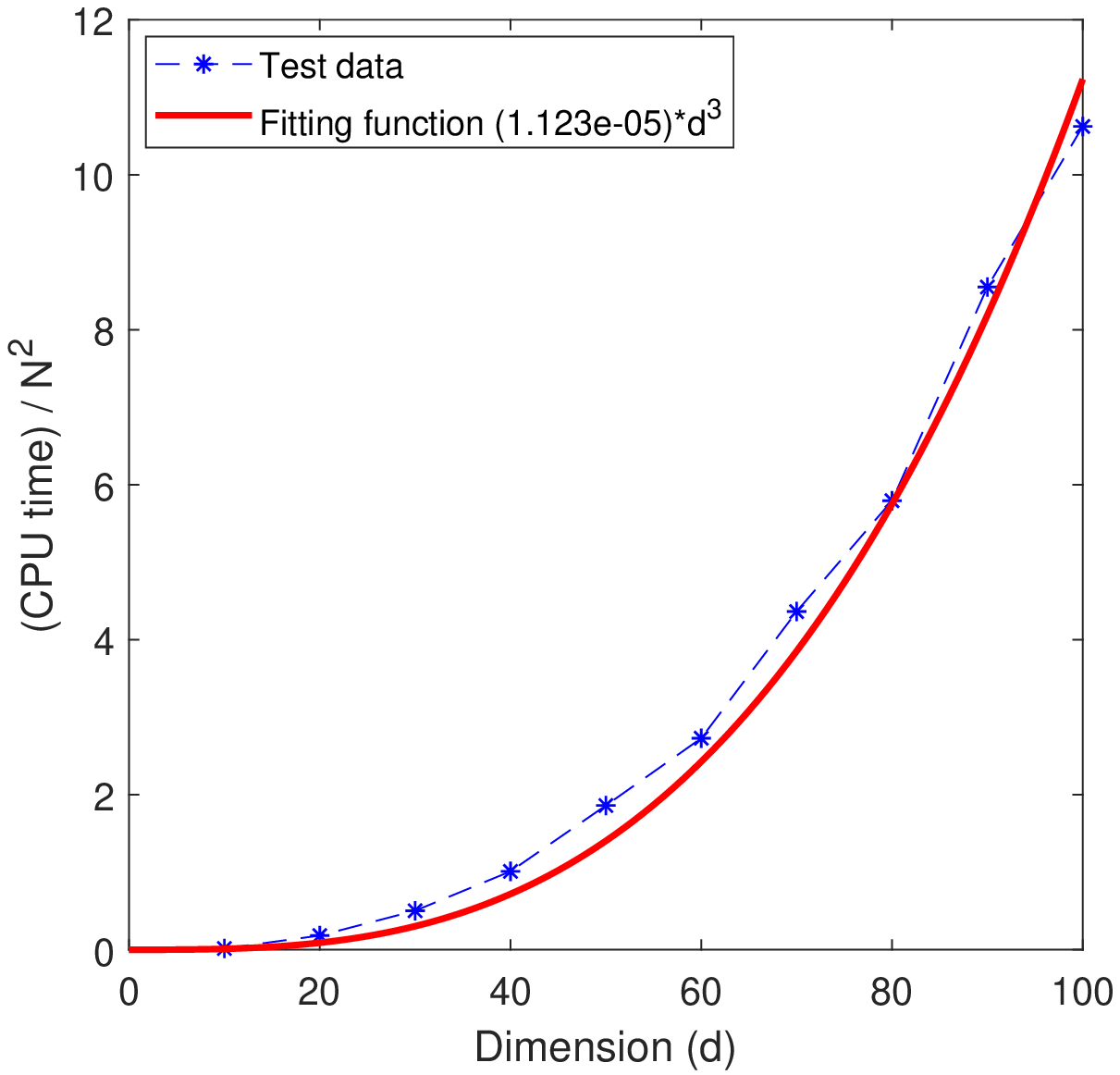}
		\includegraphics[width=1.8in,height=1.7in]{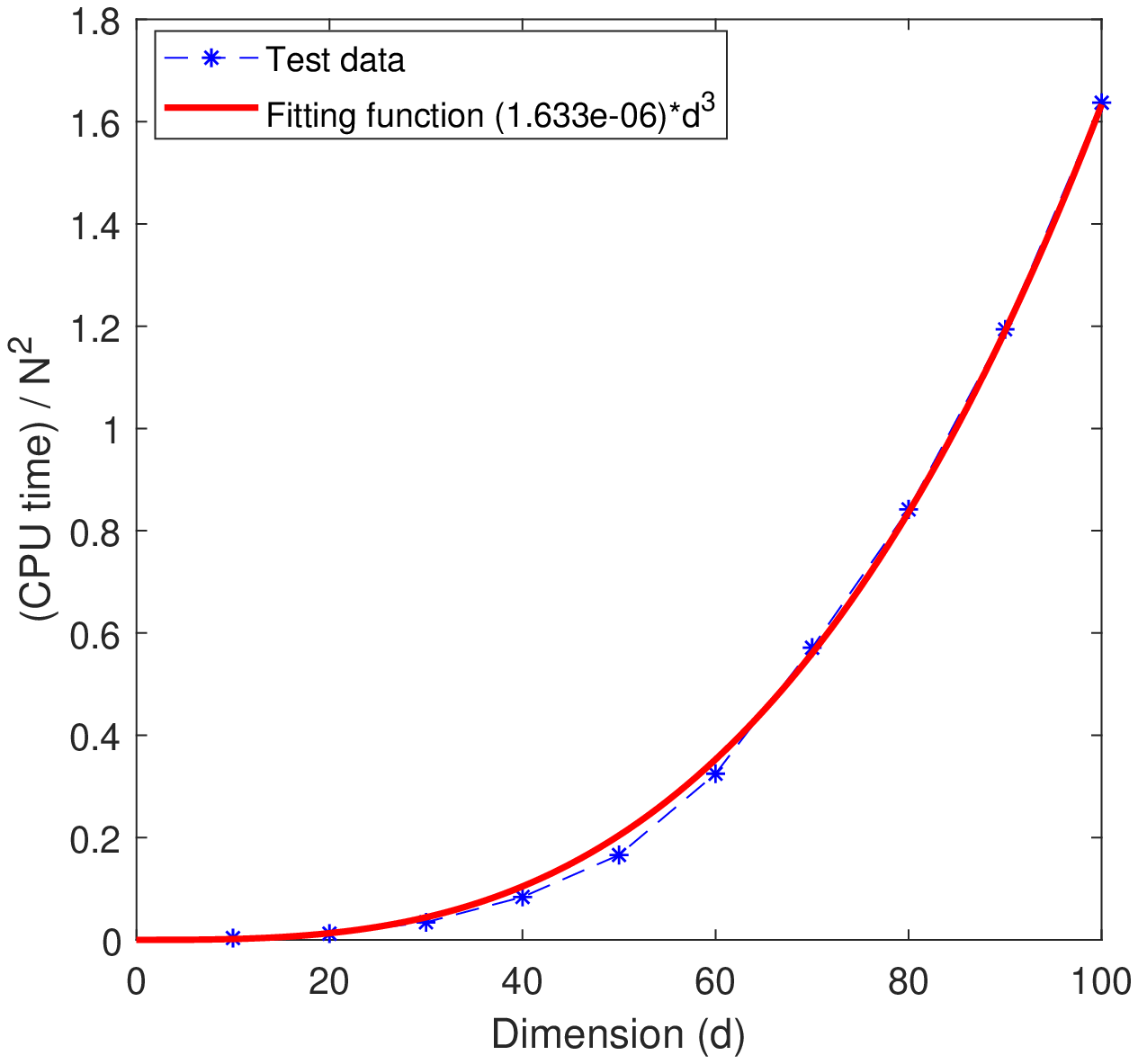}
	}
	\centerline{
		\includegraphics[width=1.8in,height=1.7in]{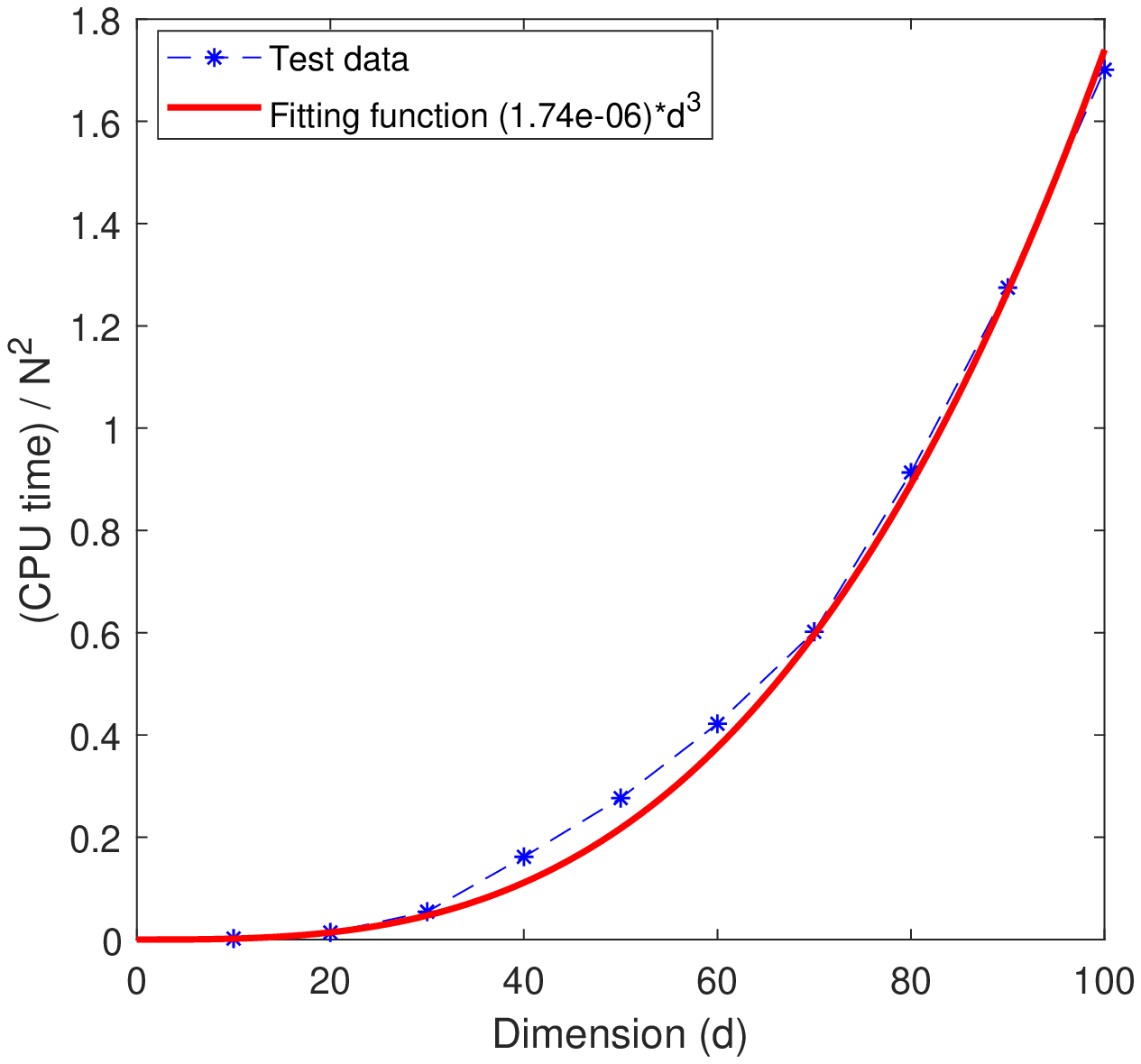}
        \includegraphics[width=1.8in,height=1.7in]{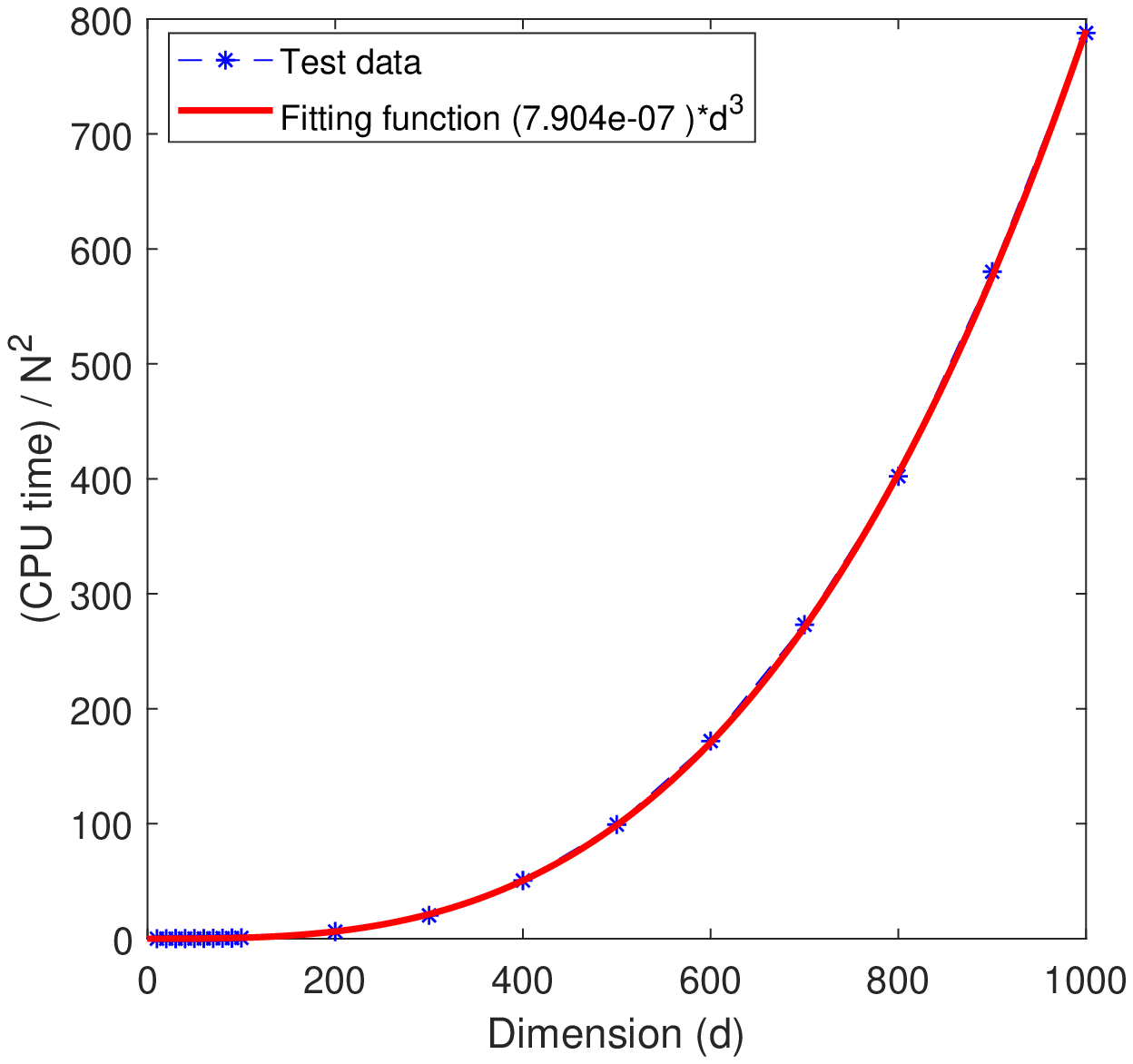}
	}  
	\centerline{
	\includegraphics[width=1.8in,height=1.7in]{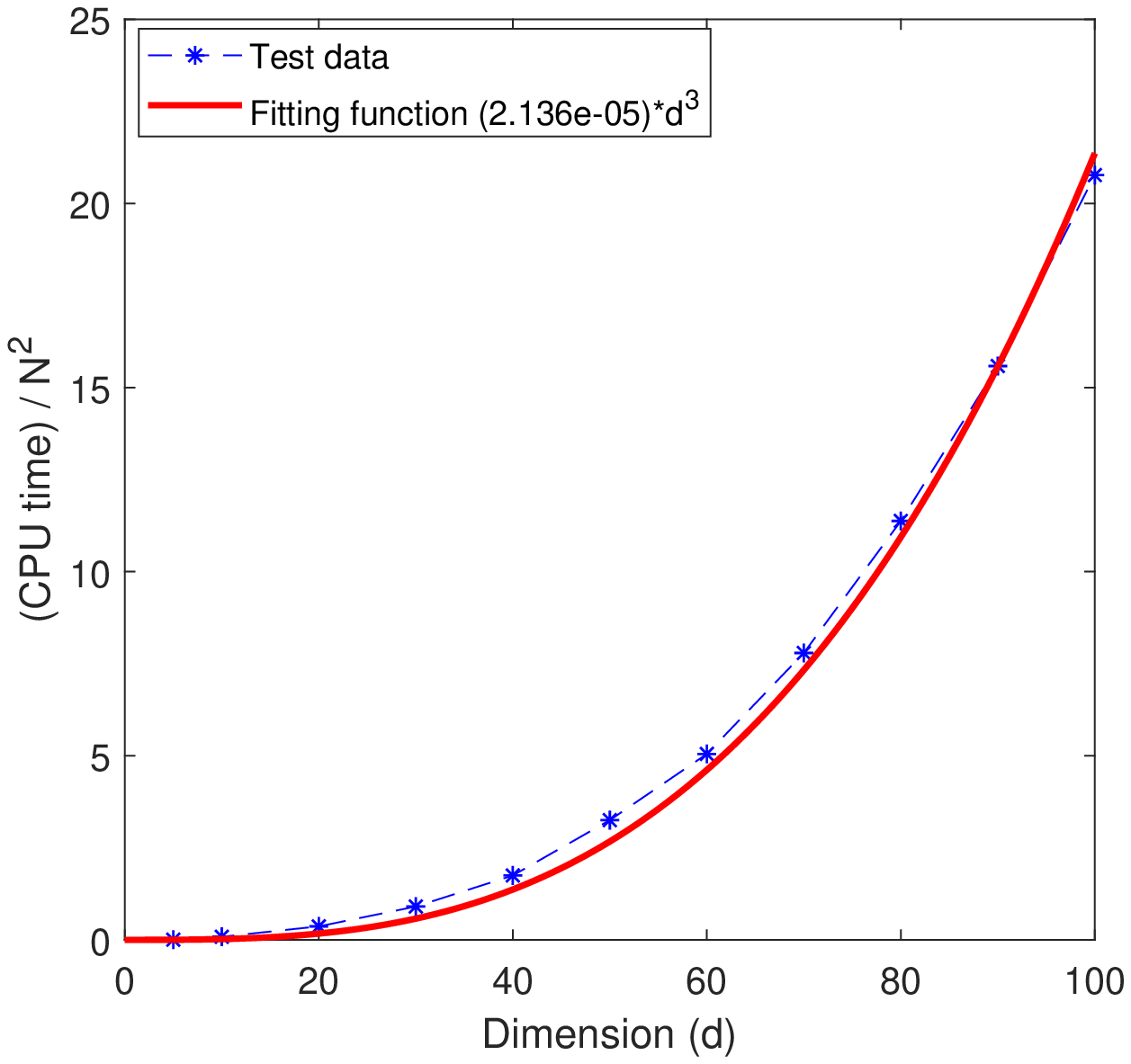}
	\includegraphics[width=1.8in,height=1.7in]{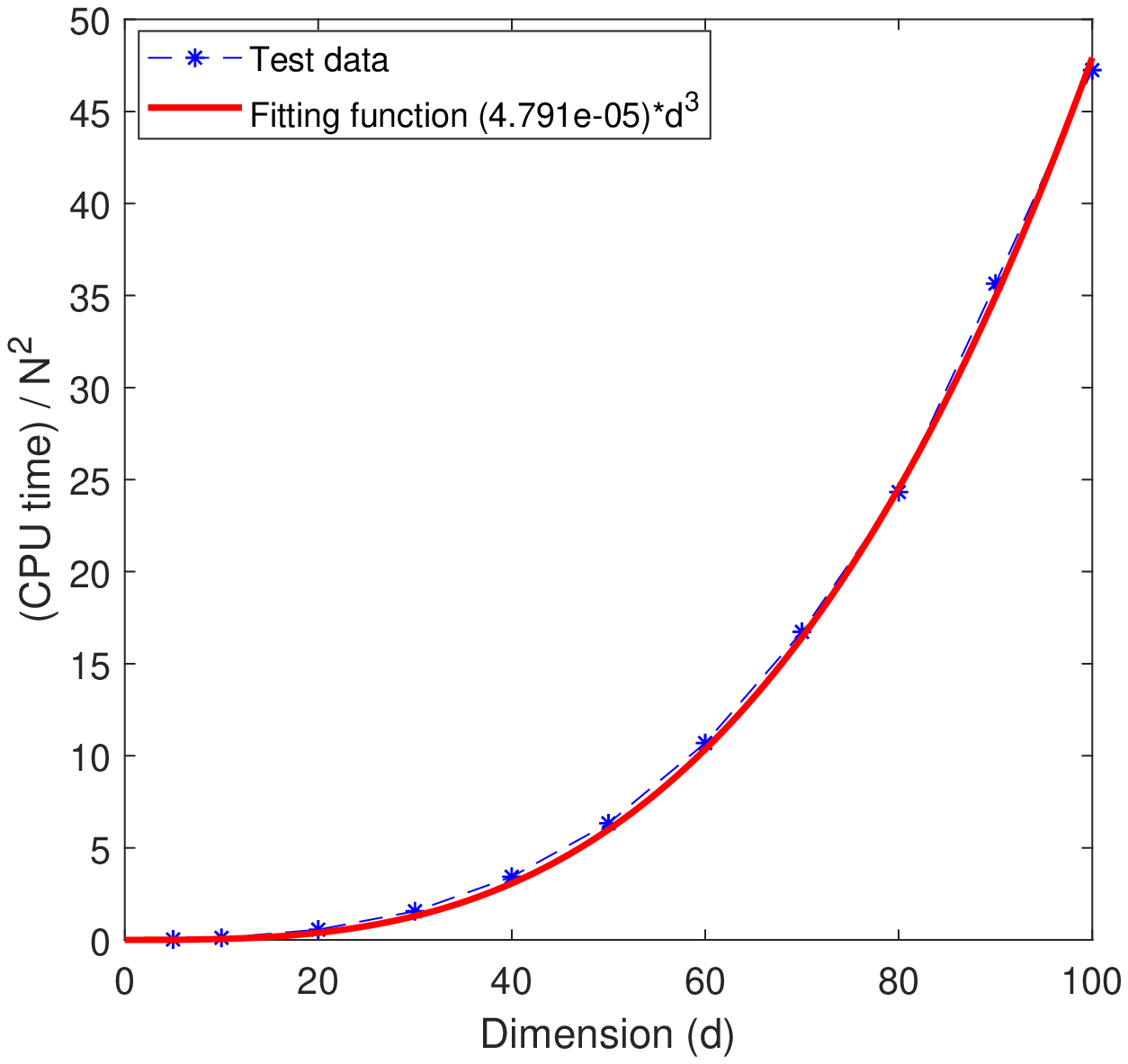}
}  
   	\caption{The relationship between the CPU time and dimension $d$.}	\label{fig.11}      
\end{figure}

\section{Conclusions}\label{sec-6}
In this paper we introduced a fast MDI (multilevel dimension iteration) algorithm (or solver) for efficiently implementing tensor product methods for high dimension numerical integration. It is based on the idea of computing the function evaluations at all integration points in cluster and iteratively along each coordinate direction, so many computations can be reused in each iteration. 
It was showed numerically based on the simulation data that the computational complexity (in terms of the CPU time) of the MDI algorithm grows at most cubically in the 
dimension $d$, and overall in the order $O(d^3N^2)$, which shows that the proposed 
MDI algorithm could effectively circumvent the curse of the dimensionality in high dimensional numerical integration, hence, makes tensor product methods not only become 
competitive but also can excel. Extensive numerical tests were provided to gauge 
the performance of the MDI algorithm and to do performance comparisons with the 
standard tensor product methods and especially with the Monte Carlo (MC) method. They
demonstrated that the MDI algorithm (regardless the choice of the 1-d base 
quadrature rules) is faster than the MC method in low and medium dimensions (i.e., $d\approx 100$), much faster in very high dimensions (i.e., $d\approx 1000$), and succeeds even when the MC method fails. 
As the idea of the MDI algorithm is applicable to any quadrature rule whose integration points have a lattice-like structure, this extension will be further investigated in the 
 future. Another direction of continuing this research is to sharpen the 
dimension-iteration idea to develop even faster algorithms which can achieve the optimal computational complexity (in terms of the CPU time) of the order $O(Nd)$, we shall present those new results in a forthcoming work in the near future. 


%









\begin{thebibliography}{99}
\bibliographystyle{abbrv}
	
\bibitem{BG14}
{\sc H.-J. Bungartz and M. Griebel},
{\em Sparse grids}, Acta Numer., 13:147--269, 2014.
	
	
\bibitem{Burden-Faires}
{\sc R. L. Burden and J. D. Faires},
{\em Numerical Analysis}, 10th edition, Cengage Learning, 2015.


\bibitem{Caflisch98}
{\sc R. E.  Caflisch}, 
{\em Monte Carlo and quasi-Monte Carlo methods}, Acta Numer., 7:1--49, 1998.

\bibitem{CDW20}
{\sc J. Chen,  R. Du, and K. Wu},
{\em A comparison study of deep Galerkin method and deep Ritz method for elliptic problems with different boundary conditions}, Commun. Math. Res., 36:354--376, 2020. 



\bibitem{DKS13}
{\sc J. Dick, F. Y. Kuo, and I. H. Sloan},
{\em High-dimensional integration: the quasi-Monte Carlo way},  Acta Numer. 22:133--288,
2013.

\bibitem{Dos12}
{\sc J. Dos Santos Azevedo and S. Pomponet Oliveira}, 
{\em  A numerical comparison between quasi-Monte Carlo and sparse grid stochastic collocation methods}, Commun. Comput. Phys., 12:1051--1069, 2012.


\bibitem{E_Yu18}
{\sc W. E and B. Yu}, 
{\em The deep Ritz method: A deep learning-based numerical algorithm for solving variational problems}, Commun. Math. and Stat., 6:1-12,  2018.


\bibitem{GH10}
{\sc M. Griebel and M. Holtz},
{\em Dimension-wise integration of high-dimensional functions with applications to finance},
J. Complexity, 26:455--489, 2010. 

\bibitem{HJE18}
{\sc J. Han, A. Jentzen, and W. E}, 
{\em Solving high-dimensional partial differential equations using deep learning}, 
PNAS, 115:8505--8510, 2018.

\bibitem{HMNR10}
{\sc F. J. Hickernell, T. M\"uller-Gronbach, B. Niu, and K. Ritter},
{\em Multi-level Monte Carlo algorithms for infinite-dimensional integration on $R^N$},  J.  Complexity  26:229--254, 2010.

\bibitem{KSS11}
{\sc F. Y. Kuo, C. Schwab, and I. H. Sloan},
{\em Quasi-Monte Carlo methods for high-dimensional integration: the standard (weighted Hilbert space) 
	setting and beyond}, ANZIAM J., 53:1--37, 2011. 

\bibitem{Lu04}
{\sc J. Lu and L. Darmofal}, 
{\em Higher-dimensional integration with Gaussian weight for applications in probabilistic design}, SIAM J. Sci. Comput., 26:613--624, 2004.

\bibitem{LMMK21}
{\sc L. Lu, X. Meng,  Z. Mao, and G. E. Karniadakis},
{\em DeepXDE: A deep learning library for solving differential equations},
SIAM Rev.,  63:208--228, 2021.

\bibitem{Ogata89}
{\sc Y. Ogata},   
{\em A Monte Carlo method for high dimensional integration},  Numer. Math. 55:137--157,  1989.

\bibitem{SS18}
{\sc J. Sirignano and K. Spiliopoulos}, 
{\em DGM: A deep learning algorithm for solving partial differential equations}, J.  Comput. Phys.
375: 1339--1364, 2018. 

\bibitem{Stoer-Bulirsch}
{\sc J. Stoer and R. Bulirsch},
{\em Introduction to Numerical Analysis},  Springer, New York,  1980.

\bibitem{Wipf13} 
{\sc A. Wipf}, 
{\em High-Dimensional Integrals}, in {\em Statistical Approach to Quantum Field Theory}, Lecture Notes in Physics, 100:25--46, Springer, 2013.

\bibitem{Xu20}
{\sc J. Xu},
{\em Finite neuron method and convergence analysis},  Commun. Comput.  Phys,  28:1707--1745, 2020.


\end{thebibliography}
 \end{document}